# A GENERALIZED DEFINITION OF CAPUTO DERIVATIVES AND ITS APPLICATION TO FRACTIONAL ODES

LEI LI* AND JIAN-GUO LIU†

**Abstract.** We propose a generalized definition of Caputo derivatives from $t = 0$ of order $\gamma \in (0,1)$ using a convolution group, and we build a convenient framework for studying initial value problems of general nonlinear time fractional differential equations. Our strategy is to define a modified Riemann–Liouville fractional calculus, which agrees with the traditional Riemann–Liouville definition for $t > 0$ but includes some singularities at $t = 0$ so that the group property holds. Then, making use of this fractional calculus, we introduce the generalized definition of Caputo derivatives. The new definition is consistent with various definitions in the literature while revealing the underlying group structure. The underlying group property makes many properties of Caputo derivatives natural. In particular, it allows us to deconvolve the fractional differential equations to integral equations with completely monotone kernels, which then enables us to prove the general comparison principle with the most general conditions. This then allows for *a priori* energy estimates of fractional PDEs. Since the new definition is valid for locally integrable functions that can blow up in finite time, it provides a framework for solutions to fractional ODEs and fractional PDEs. Many fundamental results for fractional ODEs are revisited within this framework under very weak conditions.

**1. Introduction.** The fractional calculus in continuous time has been used widely in physics and engineering for memory effect, viscoelasticity, porous media, etc. [25, 8, 18, 21, 27, 7, 1]. Some rigorous justifications of using Caputo derivatives in modeling can be found, for example, in [15, 3]. Given a function $\varphi(t)$, the fractional integral from $t = 0$ with order $\gamma > 0$ is given by Abel's formula,

$$(1) \qquad J_\gamma \varphi(t) = \frac{1}{\Gamma(\gamma)} \int_0^t \varphi(s)(t-s)^{\gamma-1} ds, \ t > 0.$$

Denote by $\theta(t)$ the Heaviside step function; then the fractional integral is just the convolution between the kernel $g_\gamma(t) = \frac{\theta(t)}{\Gamma(\gamma)} t^{\gamma-1}$ and $\theta(t)\varphi(t)$ on $\mathbb{R}$. By this fact, it is clear that the integral operators $J_\gamma$ form a semigroup. The kernel $g_\gamma$ is completely monotone if $\gamma \in (0,1]$. A function $g : (0,\infty) \to \mathbb{R}$ is completely monotone if $(-1)^n g^{(n)} \geq 0$ for $n = 0, 1, 2, \ldots$. The famous Bernstein theorem says that a function is completely monotone if and only if it is the Laplace transform of a Radon measure on $[0,\infty)$ ([28, 6]).

For the derivatives, there are two types that are commonly used: the Riemann–Liouville definition and the Caputo definition (see [25, 18]). Letting $n - 1 < \gamma < n$ where $n$ is a positive integer, the Riemann–Liouville and Caputo fractional derivatives

---

*Department of Mathematics, Duke University, Durham, NC 27708, USA. (leili@math.duke.edu).
†Departments of Physics and Mathematics, Duke University, Durham, NC 27708, USA. (jliu@phy.duke.edu).

(from $t = 0$) are defined respectively as

$$(2) \qquad D_{rl}^{\gamma}\varphi(t) = \frac{1}{\Gamma(n-\gamma)}\frac{d^n}{dt^n}\int_0^t \frac{\varphi(s)}{(t-s)^{\gamma+1-n}}ds,\ t>0,$$

$$(3) \qquad D_c^{\gamma}\varphi(t) = \frac{1}{\Gamma(n-\gamma)}\int_0^t \frac{\varphi^{(n)}(s)}{(t-s)^{\gamma+1-n}}ds,\ t>0.$$

When $\gamma = n$, the derivatives are defined to be the usual $n$th order derivatives.

Both types of fractional derivatives mentioned do not have the semigroup property (i.e. $D^{\alpha}D^{\beta}$ does not always equal $D^{\alpha+\beta}$; see Lemma 3.5). Because the integrals can be given by convolution, we expect the fractional derivative to be determined by the convolution inverse and the fractional calculus may be given by the convolution group, and therefore we desire the fractional derivatives to form a semigroup too. In [11, Chap. 1, sect. 5.5] by Gel'fand and Shilov, integrals and derivatives of arbitrary order for a distribution $\varphi \in \mathscr{D}'(\mathbb{R})$ supported on $[0,\infty)$ are defined to be the convolution

$$\varphi^{(\alpha)} = \varphi * g_\alpha, \quad g_\alpha(t) := \frac{t_+^{\alpha-1}}{\Gamma(\alpha)},\ \alpha \in \mathbb{C}.$$

Here $t_+ = \max(t,0) = \theta(t)t$ and $\frac{t_+^{\alpha-1}}{\Gamma(\alpha)}$ must be understood as distributions for $\text{Re}(\alpha) \leq 0$ (see [11, Chap. 1, sect. 3.5]). (Note that we have used different notations here from those in [11] so that they are consistent with the notations in this paper.) Clearly if $\alpha > 0$, $\varphi^{(\alpha)} = J_\alpha \varphi$ is the fractional integral as in (1). When $\alpha < 0$, they call so-defined $\varphi^{(\alpha)}$ to be derivatives of $\varphi$. As we discuss below, these derivatives are the Riemann–Liouville derivatives for $t > 0$ but may be different at $t = 0$. The calculus (integrals and derivatives) $\{\varphi^{(\alpha)} : \alpha \in \mathbb{R}\}$ forms a group.

In this paper, we will find more convenient expressions for the distributions $g_\alpha = \frac{(\theta(t)t)^{\alpha-1}}{\Gamma(\alpha)}$ when $\alpha < 0$ and then extend the ideas in [11, Chap. 1, sect. 5.5] to define fractional calculus of a certain class of distributions $\mathscr{E}^T \subset \mathscr{D}'(-\infty, T),\ T \in (0, \infty]$ (equation (13)):

$$\varphi \mapsto I_\alpha \varphi,\ \forall \varphi \in \mathscr{E}^T,\ \alpha \in \mathbb{R}.$$

The alternative expressions for $g_\alpha$ enable us to find the fractional calculus more easily, and motivate us to extend the Caputo derivatives later. If $T = \infty$ and $\text{supp}\,\varphi \subset [0, \infty)$, $I_\alpha\varphi$ is just $\varphi^{(\alpha)}$ as in [11]. We then make the distribution causal (i.e. zero when $t < 0$) by, roughly, multiplying $\theta(t)$ (see equation (20) for more rigorous treatment) and then define (Definition 2.14)

$$J_\alpha \varphi := I_\alpha(\theta\varphi),\ \alpha \in \mathbb{R}.$$

We use the same notation $J_\alpha$ as in (1) because this is exactly Abel's formula for $\alpha > 0$. $J_\alpha\varphi$ agrees with the Riemann–Liouville fractional calculus for $t > 0$, but includes some singularities at $t = 0$ (see Section 2.2.2 for more details). This means that we only need to include some singularities at $t = 0$ to make Riemann–Liouville calculus a group. These singularities are expected since the causal functions have jumps at $t = 0$. Hence, we call $J_\alpha$ the modified Riemann–Liouville calculus, which then places the foundation for us to generalize the Caputo derivatives.

Although the Caputo derivatives do not have group property, they are suitable for initial value problems and share many properties with the ordinary derivative, so they see a lot of applications in engineering ([25, 18, 7]). In the traditional definition (3), one has to define the $\gamma$th order derivative ($0 < \gamma < 1$) using the first order derivative.

This is unnatural since intuitively it can be defined for functions that are 'γth' order smooth only. In [1], Allen, Caffarelli and Vasseur have introduced an alternative form of Caputo derivative based on integration by parts for $\gamma \in (0,1)$

(4)
$$D_c^\gamma \varphi(t) = \frac{1}{\Gamma(1-\gamma)} \left( \frac{\varphi(t) - \varphi(0)}{t^\gamma} + \gamma \int_0^t \frac{\varphi(t) - \varphi(s)}{(t-s)^{\gamma+1}} ds \right)$$
$$= \frac{-1}{\Gamma(-\gamma)} \int_{-\infty}^t \frac{\tilde{\varphi}(t) - \tilde{\varphi}(s)}{(t-s)^{\gamma+1}} ds,$$

where
$$\tilde{\varphi}(s) = \begin{cases} \varphi(s), & s \geq 0, \\ \varphi(0), & s < 0. \end{cases}$$

In this definition, $\varphi'$ derivative is not needed but the integral can blow up for those $t$ where the function fails to be Hölder continuous of order $\gamma$. There are many other ways to extend the definition of Caputo derivatives for broader class of functions. For example, in [18, sect. 2.4], the Caputo derivatives are defined using the Riemann–Liouville derivatives. In [14], Gorenflo, Luchko, and Yamamoto used a functional analysis approach to extend the traditional Caputo derivatives to certain Sobolev spaces. In [2], the right and left Caputo derivatives are generalized to some integrable functions on $\mathbb{R}$ through a duality formula.

In this paper, one of our main purposes is to use the convolution group $\{g_\alpha\}$ to generalize the Caputo derivatives from $t=0$ of order $\gamma \in (0,1)$ for initial value problems. The strategy here clearly applies for any starting point of memory, not just $t=0$. Our definition reveals explicitly the underlying group structure (one can do similar things for $\gamma \notin (0,1)$ without much difficulty, but we do not consider this in this paper. See the comments in Section 3.). The idea is to remove the singular term $\frac{\varphi_0 \theta(t)}{\Gamma(1-\gamma)} t^{-\gamma}$, which corresponds to the jump of the causal function ($\varphi$ is causal if $\varphi(t) = 0$ for $t < 0$), from the modified Riemann–Liouville derivative $J_{-\gamma} \varphi$ (Definition 3.4):
$$D_c^\gamma \varphi := J_{-\gamma} \varphi - \varphi_0 \frac{\theta(t)}{\Gamma(1-\gamma)} t^{\gamma-1} = J_{-\gamma}(\varphi - \varphi_0), \ \gamma \in (0,1).$$

In general, $\varphi_0$ can be any real number and this then allows for the generalized definition. If $t=0$ is a Lebesgue point on the right and $\varphi_0 = \varphi(0+)$, by removing the singularity, the memory effect is then counted from $t=0^+$ and this is probably why Caputo derivatives are suitable for initial value problems. Counting memory from $t=0^+$ is physical. For example, in the generalized Langevin equation model derived from interacting particles (Kac-Zwanzig model) ([10, 29, 19]), the memory is from $t=0$ (A simple derivation can be found in [19, Section 4.3]). If the fractional Gaussian noise is considered, one will then have power law kernel for the memory [19], which will possibly lead to fractional calculus.

Definition 3.4 unifies all the current existing definitions, providing a new framework to study of solutions of fractional ODEs and PDEs in very weak sense and the group structure behind makes the properties of Caputo derivatives natural. For example, the definition covers the traditional definition (3) (see Proposition 3.6) and (4) used in [1] (see Corollary 3.12). Our definition shares similarity in spirit with the one used in [18], but ours is based on the modified Riemann–Liouville operators and recovers the group structure. The definition in [2] defines the Caputo derivatives through a duality relation which is valid for $L^1_{\text{loc}}(\mathbb{R})$ functions satisfying

$$\int_{\mathbb{R}} \frac{|\varphi(t)|}{1+|t|^{1+\gamma}} dt < \infty.$$

Indeed, the definition defines Caputo derivative from $t = -\infty$. If $\varphi(t) = \varphi(0)$ for all $t \leq 0$, with some work, one can show that our definition agrees with theirs (up to a negative sign). If $\varphi$ is causal ($\varphi(t) = 0$ for $t < 0$), the Caputo derivative in [2] is indeed the Riemann–Liouville derivative from $t = 0$. However, since our definition counts memory from $t = 0^+$, the values on $(-\infty, 0)$ do not matter while the derivative in [2] clearly depends on the values on $(-\infty, 0)$, which is the desired feature for initial value problems. More importantly, compared with the one in [2], Definition 3.4 is able to define Caputo derivative for functions that go to infinity in finite time, which is suitable for fractional ODEs where the solution can blow up in finite time. Another contribution of our definition is that the group structure brought by $J_\alpha$ reveals the natural properties of Caputo derivatives. For example, using $J_\alpha$, we can clearly see how $D^\alpha D^\beta$ is related to $D^{\alpha+\beta}$ (see Lemma 3.5); as another example, one can de-convolve the fractional differential equations with Caputo derivatives to integral equations with completely monotone kernel $\frac{\theta(t)}{\Gamma(\gamma)} t^{\gamma-1}$, which is the fundamental theorem of fractional calculus, without asking for extra requirements on the regularity of $\varphi(t)$ (Theorem 3.7).

Another main contribution of our work is to build a convenient framework for studying general nonlinear time fractional ODEs and PDEs. Fractional ODEs with various definitions of Caputo derivatives ((3) or the one defined through Riemann–Liouville derivatives) have already been discussed widely. See, for example, [8, 7]. With the new definition of Caputo derivatives (Definition 3.4), we are able to prove many results for the fractional ODEs under very weak conditions. For example, we are able to show the existence and uniqueness theorem in a particular subspace of locally integrable functions and prove the comparison principles with the most general conditions. We believe these results with weak conditions are useful for the analysis of fractional PDEs in the future.

The organization and main results of the paper are listed briefly as follows.

In Section 2, we extend the ideas in [11, Chap. 1, sect. 5.5] to define the modified Riemann–Liouville calculus for a certain class of distributions in $\mathscr{D}'(-\infty, T)$. Then, we prove that the modified Riemann–Liouville calculus indeed changes the regularity of functions as we expect: $J_\alpha : \tilde{H}^s(0, T) \to \tilde{H}^{s+\alpha}(0, T)$ ($s \geq \max(0, -\alpha)$) (Roughly speaking, $\tilde{H}^s(0, T)$ is the Sobolev space containing functions that are zero at $t = 0$, and one can refer to Section 2.4 for detailed explanation). In Section 3, an extension of Caputo derivatives is proposed so that both the first derivative and Hölder regularity of the function are not needed in the definition. Some properties of the new Caputo derivatives are proved, which may be used for fractional ODEs and fractional PDEs. Particularly the fundamental theorem of the fractional calculus is valid as long as the function is right continuous at $t = 0$, which allows us to transform the differential equations with orders in $(0, 1)$ to integral equations with completely monotone kernels (Theorem 3.7). In Section 4, based on the definitions and properties in Section 3, we prove some fundamental results of fractional ODEs with quite general conditions. Especially, we show the existence and uniqueness of the fractional ODEs using the fundamental theorem (Theorem 4.4), and also show the comparison principle (Theorem 4.10). This then provides the essential tools for *a priori* energy estimates of fractional PDEs. Then, as an interesting example of fractional ODEs, we show that the fractional Hamiltonian system does not preserve the 'energy', which can decay to zero.

**2. Time-continuous groups and fractional calculus.** In this section, we generalize the idea in [11] to define the fractional calculus for a particular class of distributions in $\mathscr{D}'(-\infty, T)$ and then define the modified Riemann–Liouville calculus. We first write out the time-continuous convolution group in [11, Chap. 1] in more convenient forms and then show in detail how to define the fractional calculus using convolutions. The construction in this section then provides essential tools for us to propose a generalized definition of Caputo derivatives in Section 3.

**2.1. A time-continuous convolution group.** In [11, Chap. 1], it is mentioned that $\{\frac{t_+^{\alpha-1}}{\Gamma(\alpha)} : \alpha \in \mathbb{C}\}$ can be used to define integrals and derivatives for distributions supported on $[0, \infty)$. If $\alpha \in \mathbb{R}$, these distributions form a convolution group where the convolution is well defined since all the distributions have supports that are one-side bounded. The expression $g_\alpha = \frac{t_+^{\alpha-1}}{\Gamma(\alpha)}$ with $\alpha < 0$ is not convenient for us to use. The alternative expressions for $g_\alpha$ enable us to find the fractional calculus more easily, and motivate us to extend the Caputo derivatives later.

Recall that $\theta$ represents the Heaviside step function. We note that

$$(5) \qquad \mathscr{C}_+ := \left\{ g_\alpha : g_\alpha(t) = \frac{\theta(t) t^{\alpha-1}}{\Gamma(\alpha)}, \alpha \in \mathbb{R}, \ \alpha > 0 \right\}$$

forms a semigroup of convolution. This is because

$$g_\alpha * g_\beta(t) = \frac{1}{\Gamma(\alpha)\Gamma(\beta)} \int_0^t s^{\alpha-1}(t-s)^{\beta-1} ds = \frac{B(\alpha, \beta)}{\Gamma(\alpha)\Gamma(\beta)} t^{\alpha+\beta-1} = g_{\alpha+\beta}(t),$$

where $B(\cdot, \cdot)$ is the Beta function. We now proceed to finding the convenient expressions for the convolution group generated by $\mathscr{C}_+$. For this purpose, we need the convolution between two distributions. By the theory in [11, Chap. 1], we can convolve two distributions as long as their supports are one-side bounded. This motivates us to introduce the following set of distributions

$$(6) \qquad \mathscr{E} := \{v \in \mathscr{D}'(\mathbb{R}) : \exists M_v \in \mathbb{R}, \text{supp}(v) \subset [-M_v, +\infty)\}.$$

$\mathscr{D}'(\mathbb{R})$ is the space of distributions, which is the dual of $\mathscr{D}(\mathbb{R}) = C_c^\infty(\mathbb{R})$. Clearly, $\mathscr{E}$ is a linear vector space.

The convolution between two distributions in $\mathscr{E}$ can be defined since their supports are bounded from the left. One can refer to [11] for detailed definition of the convolution. For the convenience, we mention an equivalent strategy here: pick a partition of unity for $\mathbb{R}$, $\{\phi_i : i \in \mathbb{Z}\}$ (i.e. $\phi_i \in C_c^\infty$; $0 \le \phi_i \le 1$; on any compact set $K$, there are only finitely many $\phi_i$'s that are nonzero; $\sum_i \phi_i = 1$ for all $x \in \mathbb{R}$.) The convolution can then be defined as follows.

DEFINITION 2.1. *Given $f, g \in \mathscr{E}$, we define*

$$(7) \qquad \langle f * g, \varphi \rangle := \sum_{i \in \mathbb{Z}} \langle f * (\phi_i g), \varphi \rangle, \quad \forall \varphi \in \mathscr{D} = C_c^\infty,$$

*where $\{\phi_i\}_{i=-\infty}^\infty$ is a partition of unity for $\mathbb{R}$ and $f * (\phi_i g)$ is given by the usual definition between two distributions when one of them is compactly supported.*

The following lemma is well known:

LEMMA 2.2. *(i). The definition is independent of $\{\phi_i\}$ and agrees with the usual definition of convolution between distributions whenever one of the two distributions is compactly supported.*

*(ii). For $f, g \in \mathscr{E}$, $f * g \in \mathscr{E}$.*

*(iii).*

$$f * g = g * f, \tag{8}$$

$$f * (g * h) = (f * g) * h. \tag{9}$$

*(iv). We use $D$ to mean the distributional derivative. Then, letting $f, g \in \mathscr{E}$, we have*

$$(Df) * g = D(f * g) = f * Dg. \tag{10}$$

It is well known that

LEMMA 2.3. *$g_0 = \delta(t)$ is the convolution identity and for $n \in \mathbb{N}$, $g_{-n} = D^n \delta$ is the convolution inverse of $g_n$.*

For $0 < \gamma < 1$, inspired by the fact $\mathcal{L}(g_\gamma) \sim 1/s^\gamma$ where $\mathcal{L}$ means the Laplace transform, we guess $\mathcal{L}(g_{-\gamma}) \sim s^\gamma$. Hence, the convolution inverse is guessed as $\sim D(\theta(t)t^{-\gamma})$, where $D$ is the distributional derivative. Actually, some simple computation verifies this:

LEMMA 2.4. *Letting $0 < \gamma < 1$, the convolution inverse of $g_\gamma$ is given by*

$$g_{-\gamma}(t) := \frac{1}{\Gamma(1-\gamma)} D\left(\theta(t)t^{-\gamma}\right). \tag{11}$$

*Proof.* We pick $\varphi \in C_c^\infty(\mathbb{R})$ and apply Lemma 2.2:

$$\langle D(\theta(t)t^{-\gamma}) * [\theta(t)t^{\gamma-1}], \varphi \rangle = -\langle \theta(t)t^{-\gamma} * \theta(t)t^{\gamma-1}, D\varphi \rangle$$
$$= -\langle B(1-\gamma, \gamma)\theta(t), \varphi' \rangle = -B(1-\gamma, \gamma) \int_0^\infty \varphi'(t) dt = B(1-\gamma, \gamma)\varphi(0).$$

This computation verifies that the claim is true. □

For $n < \gamma < n+1$, we define $g_{-\gamma} := D^n \delta * g_{n-\gamma}$. We have then written out the convolution group:

PROPOSITION 2.5. *Let $\mathscr{C} = \{g_\alpha : \alpha \in \mathbb{R}\}$. Then, $\mathscr{C} \subset \mathscr{E}$ and it is a convolution group under the convolution on $\mathscr{E}$ (Definition 2.1).*

*Proof.* Using the above facts and Lemma 2.2, we find that for any $\gamma > 0$, $g_{-\gamma}$ is the convolution inverse of $g_\gamma$. The fact that $\mathscr{C}_+$ forms a semigroup, the commutativity and associativity in Lemma 2.2 imply that $\mathscr{C}_- = \{g_{-\gamma}\}$ forms a convolution semigroup as well.

The group property can then be verified using the semigroup property and the fact that $g_\gamma * g_{-\gamma} = \delta$. □

**2.2. Time-continuous fractional calculus.** In this section, we use the group $\mathscr{C}$ to define the fractional calculus for a certain class of distributions in $\mathscr{D}'(-\infty, T)$ and the (modified) Riemann–Liouville fractional calculus.

**2.2.1. Fractional calculus for distributions.** In [11], the fractional derivatives for distributions supported in $[0, \infty)$ are given by $g_\alpha * \phi$. This can be easily generalized to distributions in $\mathscr{E}$:

DEFINITION 2.6. *For $\phi \in \mathscr{E}$, the fractional calculus of $\phi$, $I_\alpha : \mathscr{E} \to \mathscr{E}$ ($\alpha \in \mathbb{R}$), is given by*

$$I_\alpha \phi := g_\alpha * \phi. \tag{12}$$

REMARK 1. *It is well known that one may convolve the group $\mathscr{C}$ with $\phi \notin \mathscr{E}$ but some properties mentioned may be invalid. For example, $\phi = 1 \notin \mathscr{E}$, $g_1 = \theta$, $g_{-1} = D\delta$. Both $(\theta * D\delta) * 1$ and $\theta * (D\delta * 1)$ are defined where $\theta(t)$ is the Heaviside function, but they are not equal. The associativity is not valid.*

By Definition 2.6, it is clear that

LEMMA 2.7. *The operators $\{I_\alpha : \alpha \in \mathbb{R}\}$ form a group, and $I_{-n}\phi = D^n \phi$ ($n = 1, 2, 3, \ldots$) where $D$ is the distributional derivative.*

To get a better idea of what $I_\alpha$ is, we now pick $\phi \in C_c^\infty(\mathbb{R})$. If $\alpha \in \mathbb{Z}$, $I_\alpha$ gives the usual integral (where the integral is from $-\infty$) or derivative. For example,

$$I_1 \phi(t) = \theta * \phi(t) = \int_{-\infty}^t \phi(s) ds,$$

$$I_{-1} \phi = (D\delta) * \phi = \delta * D\phi = \phi'.$$

For $\alpha = -\gamma, 0 < \gamma < 1$ and $\phi \in C_c^\infty$, we have

$$I_{-\gamma} \phi(t) = \frac{1}{\Gamma(1-\gamma)} \frac{d}{dt} \int_{-\infty}^t \frac{\phi(s)}{(t-s)^\gamma} ds = \frac{1}{\Gamma(1-\gamma)} \int_{-\infty}^t \frac{\phi'(s)}{(t-s)^\gamma} ds.$$

Hence, $I_{-\gamma}$ gives the Riemann–Liouville derivatives from $t = -\infty$.

In many applications, the functions we study may not be defined beyond a certain time $T > 0$. For example, we may study fractional ODEs where the solution blows up at $t = T$. This requires us to define the fractional calculus for some distributions in $\mathscr{D}'(-\infty, T)$. Hence, we introduce the following set

$$\mathscr{E}^T := \{v \in \mathscr{D}'(-\infty, T) : \exists M_v \in (-\infty, T), \operatorname{supp}(v) \subset [M_v, T)\}. \tag{13}$$

$\mathscr{E}^T$ is not closed under the convolution (that means if we pick two distributions form $\mathscr{E}^T$, the convolution then is no longer in $\mathscr{E}^T$). Hence, we cannot define the fractional calculus directly as we do for $\mathscr{D}'(\mathbb{R})$. Our strategy is to push the distributions into $\mathscr{E}$ first and then pull it back.

Let $\{\chi_n\} \subset C_c^\infty(-\infty, T)$ be a sequence satisfying (i). $0 \le \chi_n \le 1$. (ii). $\chi_n = 1$ on $[-n, T - \frac{1}{n})$ (if $T = \infty$, $T - 1/n$ is defined as $\infty$). We introduce the extension operator $K_n^T : \mathscr{E}^T \to \mathscr{E}^\infty$ given by

$$\langle K_n^T v, \varphi \rangle = \langle \chi_n v, \varphi \rangle = \langle v, \chi_n \varphi \rangle, \; \forall \varphi \in C_c^\infty(\mathbb{R}), \tag{14}$$

where the last pairing is the one between $\mathscr{D}'(-\infty, T)$ and $C_c^\infty(-\infty, T)$. Denote $R^T : \mathscr{E}^\infty \to \mathscr{E}^T$ as the natural embedding operator. We define the fractional calculus as

DEFINITION 2.8. *For $\phi \in \mathscr{E}^T$, we define the fractional calculus of $\phi$ to be $I_\alpha^T : \mathscr{E}^T \to \mathscr{E}^T$ ($\alpha \in \mathbb{R}$):*

$$I_\alpha^T \phi = \lim_{n \to \infty} R^T (I_\alpha (K_n^T \phi)). \tag{15}$$

We check that the definition is well-given.

PROPOSITION 2.9. *Fix $\phi \in \mathscr{E}^T$.*

*(i). For any sequence $\{\chi_n\}$ satisfying the conditions given and $\epsilon > 0, M > 0$, there exists $N > 0$, such that $\forall n \geq N$ and $\forall \varphi \in C_c^\infty(-\infty, T)$ with $\mathrm{supp}\,\varphi \subset [-M, T-\epsilon]$,*

$$\langle K_n^T \phi, \varphi \rangle = \langle \phi, \varphi \rangle.$$

*(ii). The limit in Definition 2.8 exists under the topology of $\mathscr{D}'(-\infty, T)$. In particular, picking any partition of unity $\{\phi_i\}$ for $\mathbb{R}$, we have*

$$\langle I_\alpha^T \phi, \varphi \rangle = \sum_i \langle \phi, (g_\alpha \phi_i)(-\cdot) * \varphi \rangle, \quad \forall \varphi \in C_c^\infty(-\infty, T),\ \alpha \in \mathbb{R},$$

*and the value on the right-hand side is independent of choice of the partition of unity $\{\phi_i\}$.*

*Proof.* The proof for (i) is standard, which we omit. For (ii), we pick $\varphi \in C_c^\infty(-\infty, T)$. Then, $\forall n > 0$,

$$\langle R^T(I_\alpha(K_n^T \phi)), \varphi \rangle = \langle g_\alpha * (K_n^T \phi)), \varphi \rangle = \sum_i \langle (\phi_i g_\alpha) * (K_n^T \phi)), \varphi \rangle.$$

There are only finitely many terms in the sum. Then, for each term,

$$\langle (\phi_i g_\alpha) * (K_n^T \phi)), \varphi \rangle = \langle K_n^T \phi, \zeta_i * \varphi \rangle$$

where $\zeta_i(t) := (\phi_i g_\alpha)(-t)$ is a distribution supported in $[-N_1, 0]$ for some $N_1 > 0$. As a result, $\zeta_i * \varphi$ is $C_c^\infty(-\infty, T)$. By (i):

$$\langle I_\alpha^T \phi, \varphi \rangle = \lim_{n \to \infty} \sum_i \langle (\phi_i g_\alpha) * (K_n^T \phi)), \varphi \rangle = \sum_i \langle \phi, \zeta_i * \varphi \rangle.$$

Note that according to Lemma 2.2, the convolution is independent of the choice of partition of unity, and hence the value on the right-hand side is also independent of $\{\phi_i\}$. □

LEMMA 2.10. *$I_\alpha^T$ ($\alpha \in \mathbb{R}$) is independent of the choice of extension operators $\{K_n\}$. For any $T_1, T_2 \in (0, \infty]$ and $T_1 < T_2$,*

(16) $$R^{T_1} I_\alpha^{T_2} \phi = I_\alpha^{T_1} R^{T_1} \phi, \quad \forall \phi \in \mathscr{E}^{T_2}.$$

*Proof.* Let $\varphi \in C_c^\infty(-\infty, T_1)$. Then, we need to show

$$\lim_{n \to \infty} \langle g_\alpha * (K_n^{T_2} \phi), \varphi \rangle = \lim_{n \to \infty} \langle g_\alpha * (K_n^{T_1} R^{T_1} \phi), \varphi \rangle$$

We use the partition of unity $\{\phi_i\}$ for $\mathbb{R}$ and the equation is reduced to

$$\lim_{n \to \infty} \sum_i \langle (\phi_i g_\alpha) * (K_n^{T_2} \phi), \varphi \rangle = \lim_{n \to \infty} \sum_i \langle (\phi_i g_\alpha) * (K_n^{T_1} R^{T_1} \phi), \varphi \rangle$$

Since there are only finite terms that are nonzero for the sum, we switch the order of limit and summation. Denote $\zeta_i(t) = (\phi_i g_\alpha)(-t)$ which is supported in $(-\infty, 0)$. Then it suffices to show

$$\lim_{n \to \infty} \langle K_n^{T_2} \phi, \zeta_i * \varphi \rangle = \lim_{n \to \infty} \langle K_n^{T_1} R^{T_1} \phi, \zeta_i * \varphi \rangle$$

By Proposition 2.9, this equality is valid. □

Lemma 2.10 verifies that if $T = \infty$, $I_\alpha^T$ agrees with $I_\alpha$ in Definition 2.6.

LEMMA 2.11. $\{I_\alpha^T : \alpha \in \mathbb{R}\}$ *forms a group.*

*Proof.* By the result in Lemma 2.10, $\forall \alpha, \beta \in \mathbb{R}$,

$$I_\alpha^T(I_\beta^T \phi) = \lim_{n\to\infty} R^T(I_\alpha(R^T(I_\beta(K_n^T \phi)))) = \lim_{n\to\infty} R^T(I_{\alpha+\beta}(K_n^T \phi))) = I_{\alpha+\beta}^T \phi. \quad \square$$

DEFINITION 2.12. *If $\phi \in \mathscr{E}^T$, we can define the convolution between $g_\alpha : \alpha \in \mathbb{R}$ and $\phi$ as:*

$$(17) \qquad g_\alpha * \phi := I_\alpha^T \phi, \ \forall \phi \in \mathscr{E}^T.$$

It is easy to check that if $\phi \in \mathscr{E}^T$ is in $L^1[M, T]$ for some $M < T$ and $\alpha > 0$, we have in the Lebesgue sense that

$$g_\alpha * \phi(t) = \int_{-\infty}^t g_\alpha(t-s)\phi(s)\, ds, \ t < T.$$

By Lemma 2.11, we have

$$(18) \qquad g_\alpha * (g_\beta * \phi) = (g_\alpha * g_\beta) * \phi = g_{\alpha+\beta} * \phi.$$

From here on, when we fix $T \in (0, \infty]$, we can drop the superindex $T$ for convenience:

$$\mathscr{E}^T \text{ is denoted by } \mathscr{E}; \quad I^T \text{ is denoted by } I.$$

**2.2.2. Modified Riemann–Liouville calculus.** Definitions in Section 2.2.1 give the fractional calculus starting from $t = -\infty$. We are more interested in fractional calculus starting from $t = 0$ because the memory is usually counted from $t = 0$ in many applications and initial value problems are considered. In other words, for $\varphi \in \mathscr{E}$ (with $T \in (0, \infty]$ fixed), we hope to define the fractional calculus from $t = 0$. Inspired by [11, Chap. 1, sect. 5.5], we need to modify the distribution to be supported in $[0, T)$.

Consider causal distributions ('zero for $t < 0$'):

$$(19) \qquad \mathscr{G}_c := \{\phi \in \mathscr{E} \subset \mathscr{D}'(-\infty, T) : \operatorname{supp} \phi \subset [0, T)\}.$$

We now consider the causal correspondence for a general distribution in $\mathscr{E}$. Let $u_n \in C_c^\infty(-1/n, T)$ where $n = 1, 2, \ldots$ be a sequence satisfying (1) $0 \leq u_n \leq 1$, (2) $u_n(t) = 1$ for $t \in (-1/(2n), T - 1/(2n))$. Introduce the space

$$(20) \quad \mathscr{G} := \{\varphi \in \mathscr{E} : \exists \phi \in \mathscr{G}_c, u_n \varphi \to \phi \text{ in } \mathscr{D}'(-\infty, T) \text{ for any such sequence } \{u_n\}\}.$$

For $\varphi \in \mathscr{G}$, the corresponding distribution $\phi$ is denoted as $\theta\varphi$ where $\theta$ is the Heaviside step function. Clearly, if $\varphi \in L^1_{\text{loc}}(-\infty, T)$, where the notation $L^1_{\text{loc}}(U)$ represents the set of all locally integrable function defined on $U$, $\theta\varphi$ can be understood as the usual multiplication.

LEMMA 2.13. $\mathscr{G}_c \subset \mathscr{G}$. $\forall \varphi \in \mathscr{G}_c, \theta\varphi = \varphi$.

This claim is easy to verify and we omit the proof. This then motivates the following definition:

DEFINITION 2.14. *The (modified) Riemann–Liouville operators $J_\alpha : \mathcal{G} \to \mathcal{G}_c$ are given by*

$$J_\alpha \varphi := I_\alpha(\theta\varphi) = g_\alpha * (\theta\varphi), \tag{21}$$

*where $g_\alpha * (\theta\varphi)$ is understood as in Equation* (17).

PROPOSITION 2.15. *Fix $\varphi \in \mathcal{E}$. $\forall \alpha, \beta \in \mathbb{R}$, $J_\alpha J_\beta \varphi = J_{\alpha+\beta}\varphi$ and $J_0\varphi = \theta\varphi$. If we make the domain of them to be $\mathcal{G}_c$ (i.e. the set of causal distributions), then they form a group.*

*Proof.* One can verify that $\mathrm{supp}(J_\alpha\varphi) \subset [0, T)$. Hence, $\theta J_\alpha\varphi = J_\alpha\varphi$. The claims follow from the properties of $I_\alpha$. If $\varphi \in \mathcal{G}_c$, then $\varphi$ is identified with $\theta(t)\varphi$. □

We are more interested in the cases where $\varphi$ is locally integrable. We call them modified Riemann–Liouville because for good enough $\varphi$ they agree with the traditional Riemann–Liouville operators (Equation (2)) at $t > 0$ while there are some extra singularities at $t = 0$. Now, let us illustrate this by checking some special cases.

When $\alpha > 0$ and $\varphi$ is a continuous function, we have verified that (21) gives the Abel's formula of fractional integrals (1). It would be interesting to look at the formulas for $\alpha < 0$ and smooth $\varphi$:

- When $-1 < \alpha < 0$, we have for any $t < T$

$$J_\alpha\varphi(t) = \frac{\theta(t)}{\Gamma(1-\gamma)} D \int_0^t \frac{\varphi(s)}{(t-s)^\gamma} ds$$

$$= \frac{1}{\Gamma(1-\gamma)} (\theta(t)t^{-\gamma}) * (\theta(t)\varphi' + \delta(t)\varphi(0)) \tag{22}$$

$$= \frac{1}{\Gamma(1-\gamma)} \int_0^t \frac{1}{(t-s)^\gamma} \varphi'(s)\, ds + \varphi(0) \frac{\theta(t)}{\Gamma(1-\gamma)} t^{-\gamma},$$

where $\gamma = -\alpha$. This is the Riemann–Liouville fractional derivative.

- When $\alpha = -1$, we have

$$J_{-1}\varphi(t) = D\Big(\theta(t)\varphi(t)\Big) = \theta(t)\varphi'(t) + \delta(t)\varphi(0). \tag{23}$$

We can verify easily that $J_{-1}J_1\varphi = J_1J_{-1}\varphi = \varphi$. Traditionally, the Riemann–Liouville derivatives for integer values are defined as the usual derivatives. $J_{-1}\varphi(t)$ agrees with the usual derivative for $t > 0$ but it has a singularity due to the jump of $\theta(t)\varphi$ at $t = 0$.

- When $\alpha = -1 - \gamma$. By the group property, we have for $t < T$

$$J_\alpha\varphi(t) = J_{-1}(J_{-\gamma}\varphi) = \frac{1}{\Gamma(1-\gamma)} D\Big(\theta(t) D \int_0^t \frac{\varphi(s)}{(t-s)^\gamma}\, ds\Big)$$

$$= \frac{1}{\Gamma(2-|\alpha|)} D\Big(\theta(t) D \int_0^t \frac{1}{(t-s)^{|\alpha|-1}} \varphi(s)\, ds\Big).$$

This is again the Riemann–Liouville derivative for $t > 0$.

We then call $\{J_\alpha\}$ the **modified Riemann–Liouville operators**. Clearly, $J_\alpha\varphi$ agrees with the traditional Riemann–Liouville calculus as distributions in $\mathscr{D}'(0,\infty)$ (i.e. they agree for $t > 0$). However, at $t = 0$, there is some difference. For example, $J_{-1}$ gives an atom $\varphi(0)\delta(t)$ at the origin so that $J_1J_{-1} = J_{-1}J_1 = J_0$. The singularities at $t = 0$ are expected since the causal function $\theta(t)\varphi(t)$ usually has a jump at $t = 0$. (If we use the Riemann–Liouville calculus from $t = -\infty$, it forms a group automatically and Riemann–Liouville derivative with order $n$ ($n \in \mathbb{N}_+$) agrees with the usual $n$th order derivative without singular terms.)

### 2.3. Another group for right derivatives.
Now consider another group $\widetilde{\mathscr{C}}$ generated by

$$\tilde{g}_\alpha(t) := \frac{\theta(-t)}{\Gamma(\alpha)}(-t)^{\alpha-1}, \quad \alpha > 0. \tag{24}$$

For $0 < \gamma < 1$, $\tilde{g}_{-\gamma} = -\frac{1}{\Gamma(1-\gamma)} D(\theta(-t)(-t)^{-\gamma})$ ($D$ means the derivative on $t$ and $D\theta(-t) = -\delta(t)$). The action of this group is well defined if we act it on distributions that have supports on $(-\infty, M]$ or on functions that decay faster than rational functions at $\infty$.

This group can generate fractional derivatives that are noncausal. For example if $\phi \in C_c^\infty(\mathbb{R})$,

$$(\tilde{g}_{-\gamma} * \phi)(t) = -\frac{1}{\Gamma(1-\gamma)} \frac{d}{dt} \int_t^\infty (s-t)^{-\gamma} \phi(s) \, ds. \tag{25}$$

This derivative is called the **right Riemann–Liouville derivative** in some literature (See e.g. [18]). The derivative at $t$ depends on the values in the future and it is therefore noncausal.

This group is actually the dual of $\mathscr{C}$ in the following sense

$$\langle g_\alpha * \phi, \varphi \rangle = \langle \phi, \tilde{g}_\alpha * \varphi \rangle, \tag{26}$$

where both $\phi$ and $\varphi$ are in $C_c^\infty(\mathbb{R})$. (If $\phi$ and $\varphi$ are not compactly supported or do not decay at infinity, then at least one group is not well defined for them.) This dual identity actually provides a type of integration by parts.

It is interesting to write explicitly out the case $\alpha = -\gamma$.

$$\int_{-\infty}^\infty g_\alpha * \phi(t) \varphi(t) \, dt = \int_{-\infty}^\infty \frac{1}{\Gamma(1-\gamma)} \int_{-\infty}^t (t-s)^{-\gamma} D\phi(s) \, ds \, \varphi(t) \, dt$$

$$= -\int_{-\infty}^\infty \frac{1}{\Gamma(1-\gamma)} \frac{D}{Ds} \int_s^\infty (t-s)^{-\gamma} \varphi(t) \, dt \, \phi(s) \, ds = \int_{-\infty}^\infty \tilde{g}_\alpha * \varphi(s) \phi(s) \, ds.$$

REMARK 2. *Alternatively, one may define the operator $I_\alpha$ by $\langle I_\alpha \phi, \varphi \rangle = \langle \phi, \tilde{g}_\alpha * \varphi \rangle$ for $\phi \in \mathscr{D}'(\mathbb{R})$ and $\varphi \in \mathscr{D}(\mathbb{R})$ whenever this is well defined. This definition however is also generally only valid for $\phi \in \mathscr{E}$. This is because $\tilde{g}_\alpha * \varphi$ is supported on $(-\infty, M]$ for some $M$. If $\phi \notin \mathscr{E}$, the definition may not make sense.*

### 2.4. Regularities of the modified Riemann–Liouville operators.
By the definition, it is expected that $\{J_\alpha\}$ indeed improve or reduce regularities as the ordinary integrals or derivatives do. In this section, we check this topic by considering their actions on a specific class of Sobolev spaces.

Let us fix $T \in (0, \infty)$ (we are not considering $T = \infty$). Recall that $H_0^s(0, T)$ is the closure of $C_c^\infty(0, T)$ under the norm of $H^s(0, T)$ ($H^s(0, T)$ itself equals the closure of $C^\infty[0, T]$). We would like to avoid the singularities that may appear at $t = 0$ but we do not require much at $t = T$. We therefore introduce the space $\tilde{H}^s(0, T)$ which is the closure of $C_c^\infty(0, T]$ under the norm of $H^s(0, T)$ (If $\varphi \in C_c^\infty(0, T]$, $\operatorname{supp} \varphi \subset C(0, T]$ and $\varphi \in C^\infty[0, T]$. $\varphi(T)$ may be nonzero.).

We now introduce some lemmas for our further discussion:

LEMMA 2.16. *Let $s \in \mathbb{R}, s \geq 0$.*

- The restriction mapping is bounded from $H^s(\mathbb{R})$ to $H^s(0,T)$, i.e. $\forall v \in H^s(\mathbb{R})$, then $v \in H^s(0,T)$ and there exists $C = C(s,T)$ such that $\|v\|_{H^s(0,T)} \leq \|v\|_{H^s(\mathbb{R})}$.
- For $v \in \tilde{H}^s(0,T)$, $\exists v_n \in C_c^\infty(\mathbb{R})$ such that the following conditions hold: (i). $\operatorname{supp} v_n \subset (0, 2T)$. (ii). $\|v_n\|_{H^s(\mathbb{R})} \leq C\|v_n\|_{H^s(0,T)}$, where $C = C(s,T)$. (iii). $v_n \to v$ in $H^s(0,T)$.

We use $\mathscr{D}'(0,T)$ to mean the dual of $\mathscr{D}(0,T) = C_c^\infty(0,T)$. Recall Definition 2.14.

LEMMA 2.17. *If $v_n \to f$ in $H^s(0,T)$, ($s \geq 0$), then, $J_\alpha v_n \to J_\alpha f$ in $\mathscr{D}'(0,T)$, $\forall \alpha \in \mathbb{R}$.*

*Proof.* Since $v_n, f$ are in $H^s$, then they are locally integrable functions.
Let $\varphi \in \mathscr{D}(0,T)$. Let $\{\phi_i\}$ a partition of unity for $\mathbb{R}$. By the definition of $J_\alpha$,

$$\langle g_\alpha * (\theta(t)v_n), \varphi \rangle = \sum_i \langle (g_\alpha \phi_i) * (\theta(t)v_n), \varphi \rangle = \sum_i \langle \theta(t)v_n, h_\alpha^i * \varphi \rangle,$$

where $h_\alpha^i(t) = (g_\alpha \phi_i)(-t)$. Note that there are only finitely many terms that are nonzero in the sum since $\varphi$ is compactly supported and $g_\alpha$ is supported in $[0, \infty)$. Since the support of $g_\alpha \phi_i$ is in $[0, \infty)$, then the support of $h_\alpha^i * \varphi$ is in $(-\infty, T)$. Further, $h_\alpha^i * \varphi \in C_c^\infty(\mathbb{R})$. Hence, in the distributional sense,

$$\langle \theta(t)v_n, h_\alpha^i * \varphi \rangle = \int_0^T v_n(t)(h_\alpha^i * \varphi)(t)dt \to \int_0^T f(t)(h_\alpha^i * \varphi)(t)dt$$
$$= \langle \theta(t)f, h_\alpha^i * \varphi \rangle = \langle (g_\alpha \phi_i) * (\theta(t)f), \varphi \rangle.$$

This verifies the claim. □

We now consider the action of $J_\alpha$ on $\tilde{H}^s(0,T)$ and we actually have:

THEOREM 2.18. *If $\min\{s, s+\alpha\} \geq 0$, then $J_\alpha$ is bounded from $\tilde{H}^s(0,T)$ to $\tilde{H}^{s+\alpha}(0,T)$. In other words, if $f \in \tilde{H}^s(0,T)$, then $J_\alpha f \in \tilde{H}^{s+\alpha}(0,T)$ and there exists a constant $C$ depending on $T$, $s$ and $\alpha$ such that*

(27) $$\|J_\alpha f\|_{H^{s+\alpha}(0,T)} \leq C\|f\|_{H^s(0,T)}, \ \forall f \in \tilde{H}^s(0,T).$$

About this topic, some partial results can be found in [18, 17, 13].

*Proof.* In the proof here, we use $C$ to mean a generic constant, i.e. $C$ may represent different constants from line to line, but we just use the same notation.

$\alpha = 0$ is trivial as we have the identity map.

Consider $\alpha < 0$ first. For $\alpha = -n$ ($n = 1, \ldots$), let $v \in C_c^\infty(0, \infty)$. $J_{-n}v \in C_c^\infty(0, \infty)$ because in this case, the action is the usual $n$th order derivative. It is clear that

$$\|J_{-n}v\|_{H^{s-n}(\mathbb{R})} \leq C\|v\|_{H^s(\mathbb{R})}.$$

Taking a sequence $v_i \in C_c^\infty$ and $\operatorname{supp} v_i \subset (0, 2T)$ such that $\|v_i\|_{H^s(\mathbb{R})} \leq C\|v_i\|_{H^s(0,T)}$, and $v_i \to f$ in $H^s(0,T)$. It then follows that $\|J_{-n}v_i\|_{H^{s-n}(\mathbb{R})} \leq C\|v_i\|_{H^s(0,T)}$. Since the restriction is bounded from $H^{s-n}(\mathbb{R})$ to $H^{s-n}(0,T)$, $J_{-n}v_i$ is a Cauchy sequence in $\tilde{H}^{s-n}(0,T)$. The limit in $\tilde{H}^{s-n}(0,T)$ must be $J_{-n}f$ by Lemma 2.17. Hence, $J_{-n}$ sends $\tilde{H}^s(0,T)$ to $\tilde{H}^{s-n}(0,T)$.

By the group property, it suffices to consider $-1 < \alpha < 0$ for fractional derivatives. Let $\gamma = |\alpha|$. We pick first $v \in C_c^\infty(0, 2T)$. We have

$$J_{-\gamma}v(t) = \frac{1}{\Gamma(1-\gamma)} \frac{d}{dt} \int_0^t s^{-\gamma} v(t-s)\, ds = \frac{1}{\Gamma(1-\gamma)} \int_0^t (t-s)^{-\gamma} v'(s)\, ds.$$

Since $J_{-\gamma}v = \frac{1}{\Gamma(1-\gamma)}(\theta(t)t^{-\gamma}) * (v')$ and $v' \in C_c^\infty(0, 2T)$, $J_{-\gamma}v$ is $C^\infty$ and $\mathrm{supp}\, J_{-\gamma}v \subset (0, \infty)$. Note that the last term is the Caputo derivative. The Caputo derivative equals the Riemann–Liouville derivative for $v \in C_c^\infty(0, 2T)$.

Since $|\mathcal{F}(\theta(t)t^{-\gamma})| \leq C|\xi|^{\gamma-1}$, we find $|\mathcal{F}(J_{-\gamma}v)| \leq C|\xi|^\gamma |\hat{v}(\xi)|$. Here, $\mathcal{F}$ represents the Fourier transform operator, while $\hat{v}$ is the Fourier transform of $v$. Hence,

$$\int (1+|\xi|^2)^{(s-\gamma)} |\mathcal{F}(J_{-\gamma}v)|^2 d\xi \leq C \int (1+|\xi|^2)^s |\hat{v}(\xi)|^2 d\xi,$$

or $\|J_{-\gamma}v\|_{H^{s-\gamma}(\mathbb{R})} \leq C\|v\|_{H^s(\mathbb{R})}$. By Lemma 2.16, the restriction is bounded

$$\|J_{-\gamma}v\|_{H^{s-\gamma}(0,T)} \leq C\|v\|_{H^s(\mathbb{R})}.$$

Take $v_i \in C_c^\infty(0, 2T)$ such that $v_i \to f$ in $H^s(0, T)$ and $\|v_i\|_{H^s(\mathbb{R})} \leq C\|v_i\|_{H^s(0,T)}$, then $\|J_{-\gamma}v_i\|_{H^{s-\gamma}(0,T)} \leq C\|v_i\|_{H^s(0,T)}$ and $J_{-\gamma}v_i$ is a Cauchy sequence in $\tilde{H}^s(0, T) \subset H^s(0, T)$. The limit in $\tilde{H}^s(0, T)$ must be $J_{-\gamma}f$ by Lemma 2.17. Hence, the claim follows for $-1 < \alpha < 0$.

Consider $\alpha > 0$ and $n \leq \alpha < n+1$. Note that $J_n$ sends $\tilde{H}^s(0, T)$ to $\tilde{H}^{s+n}(0, T)$ since this is the usual integral. We therefore only have to prove the claim for $0 < \alpha < 1$ by the group property.

For $0 < \alpha < 1$, $J_\alpha v = \frac{1}{\Gamma(\alpha)} \int_0^t s^{\alpha-1} v(t-s) ds \in C^\infty(0, \infty)$ and $\mathrm{supp}\, J_\alpha v \subset (0, \infty)$ for $v \in C_c^\infty(0, 2T)$. We again set $\gamma = |\alpha| = \alpha$. The Fourier transform of $J_\gamma v$ is $\hat{v}/(i\xi)^\gamma$. There is singularity at $\xi = 0$ because $J_\gamma v \sim t^{\gamma-1}$ as $t \to \infty$. Since we care about the behavior on $(0, T)$, we can pick a cutoff function $\zeta = \beta(x/T)$ where $\beta = 1$ on $[-1, 1]$ and zero for $|x| > 2$. $\hat{\zeta}$ is a Schwartz function.

Noting $|\mathcal{F}(\zeta J_\gamma v)| \leq |\hat{\zeta} * \hat{v}|\xi|^{-\gamma}| \leq |\hat{\zeta} * |\xi|^{-\gamma}|\|\hat{v}\|_\infty \leq C\|\hat{v}\|_\infty$, we find

$$\|\zeta J_\gamma v\|_{H^{s+\gamma}(\mathbb{R})}^2 \leq \int_\mathbb{R} (1+|\xi|^2)^{s+\gamma} |\hat{\zeta} * (\hat{v}|\xi|^{-\gamma})|^2 d\xi = \int_{|\xi|<R} + \int_{|\xi|\geq R} \leq C\|\hat{v}\|_\infty^2 + \int_{|\xi|\geq R}.$$

For $|\xi| \geq R$, we split the convolution $\hat{\zeta} * (\hat{v}|\xi|^{-\gamma})$ into two parts and apply the inequality $(a+b)^2 \leq 2(a^2 + b^2)$. It then follows that

$$\int_{|\xi|\geq R} \leq C \int_{|\xi|\geq R} d\xi (1+|\xi|^2)^{s+\gamma} \left( \left( \int_{|\eta|\geq |\xi|/2} |\hat{\zeta}(\xi-\eta)||\hat{v}(\eta)||\eta|^{-\gamma} d\eta \right)^2 \right.$$

$$\left. + \left( \int_{|\eta|\leq |\xi|/2} |\hat{\zeta}(\xi-\eta)||\hat{v}(\eta)||\eta|^{-\gamma} d\eta \right)^2 \right) = I_1 + I_2.$$

For $I_1$ by Hölder inequality and Fubini theorem,

$$I_1 \leq C \int_{|\xi|\geq R} d\xi (1+|\xi|^2)^{s+\gamma} \int_{|\eta|\geq |\xi|/2} |\hat{\zeta}(\xi-\eta)||\hat{v}(\eta)|^2 |\eta|^{-2\gamma} d\eta$$

$$\leq C \int_{|\eta|\geq R/2} d\eta |\hat{v}(\eta)|^2 |\eta|^{-2\gamma} \int_\xi |\hat{\zeta}(\xi-\eta)|(1+|\xi|^2)^{s+\gamma} d\xi$$

$$\leq C \int_{|\eta|\geq R/2} d\eta |\hat{v}(\eta)|^2 |\eta|^{-2\gamma} (1+|\eta|^2)^{\gamma+s} \leq C\|v\|_{H^s(\mathbb{R})}^2.$$

Here, $C$ depends on $R$ and $\zeta$.

For $I_2$ part, we note that $|\hat{\zeta}(\xi - \eta)| \leq C|\xi|^{-N}$ if $R$ is large enough, since $\hat{\zeta}$ is a Schwartz function.

$$\int_{|\eta| \leq |\xi|/2} |\hat{\zeta}(\xi - \eta)\hat{v}(\eta)||\eta|^{-\gamma} d\eta \leq C\|\hat{v}\|_\infty |\xi|^{-N} \int_{|\eta| \leq |\xi|/2} |\eta|^{-\gamma} d\eta \leq C\|\hat{v}\|_\infty |\xi|^{-N+1-\gamma}.$$

Hence, $I_2 \leq C\|\hat{v}\|_\infty^2$.

Overall, we have

$$\|\zeta J_\gamma v\|_{H^{s+\gamma}(\mathbb{R})} \leq C(\|\hat{v}\|_\infty + \|v\|_{H^s(\mathbb{R})}) \leq C\|v\|_{H^s(\mathbb{R})}.$$

Note that $v$ is supported in $(0, 2T)$ and $\|\hat{v}\|_\infty \leq \|v\|_{L^1(0,2T)}$, which is bounded by its $L^2(0, 2T)$ norm and thus $H^s(\mathbb{R})$ norm. The constant $C$ depends on $T$ and $\zeta$. Using again that the restriction map is bounded, we find that

$$\|J_\gamma v\|_{H^{s+\gamma}(0,T)} \leq C\|\zeta J_\gamma v\|_{H^{s+\gamma}(\mathbb{R})} \leq C\|v\|_{H^s(\mathbb{R})}.$$

The claim is true for $C_c^\infty(0, 2T)$. Again, using an approximation sequence $v_i \in C_c^\infty(0, 2T)$, $\|v_i\|_{H^s(\mathbb{R})} \leq C\|v_i\|_{H^s(0,T)}$ implies that it is true for $\tilde{H}^s(0, T)$ also. □

Enforcing $\varphi$ to be in $\tilde{H}^s(0, T)$ removes the singularities at $t = 0$. This then allows us to obtain the regularity estimates above and the Caputo derivatives will be the same as Riemann–Liouville derivatives. If $v \in \tilde{H}^0(0, T) = L^2(0, T)$, then the value of $J_\gamma v$ at $t = 0$ is well defined for $\gamma > 1/2$, which should be zero (See also [17]), because the Hölder inequality implies $\int_0^t (t - s)^{\gamma-1} v(s) ds \leq C(v) t^{\gamma-1/2}$. Actually, $\tilde{H}^\gamma(0, T) \subset C^0[0, T]$ if $\gamma > 1/2$.

**3. An extension of Caputo derivatives.** As mentioned in the introduction, we are going to propose a generalized definition of Caputo derivatives without using derivatives of integer order as in (3). Observing the calculation like (22), the Caputo derivatives $D_c^\gamma \varphi$ ($\gamma > 0$) (Equation (3)) may be defined using $J_{-\gamma}\varphi$ and the terms like $\frac{\varphi(0)}{\Gamma(1-\gamma)} t^{-\gamma}$, and hence may be generalized to a function $\varphi$ such that only $\varphi^{(m)}, m \leq [\gamma]$ exist in some sense, where $[\gamma]$ means the largest integer that does not exceed $\gamma$. We then do not need to require that $\varphi^{([\gamma]+1)}$ exist. In this paper, we only deal with $0 < \gamma < 1$ cases as they are mostly used in practice. (For general $\gamma > 1$, one has to remove singular terms related to $\varphi(0), \ldots, \varphi^{[\gamma]}(0)$, the jumps of the derivatives of $\theta(t)\varphi$, from $J_{-\gamma}\varphi$.) We prove some basic properties of the extended Caputo derivatives according to our definition, which will be used for the analysis of fractional ODEs in Section 4, and may be possibly used for fractional PDEs.

Recall that $L_{\text{loc}}^1[0, T)$ is the set of locally integrable functions on $[0, T)$, i.e. the functions are integrable on any compact set $K \subset [0, T)$. We first introduce a result from real analysis ([24]):

LEMMA 3.1. *Suppose $f, g \in L_{\text{loc}}^1[0, T)$ where $T \in (0, \infty]$, then $h(x) = \int_0^x f(x - y)g(y) dy$ is defined for almost every $x \in [0, T)$ and $h \in L_{\text{loc}}^1[0, T)$.*

*Proof.* Fix $M \in (0, T)$. Denote $\Omega = \{(x, y) : 0 \leq y \leq x \leq M\}$. $F(x, y) = |f(x - y)||g(y)|$ is measurable and nonnegative on $\Omega$. Tonelli's theorem ([24, sect. 12.4]) indicates that

$$\iint_D F(x, y) dA = \int_0^M |g(y)| \int_y^M |f(x - y)| dx dy \leq C(M),$$

for some $C(M) \in (0, \infty)$. This means that $F(x,y)$ is integrable on $D$. Hence, $h(x)$ is defined for almost every $x \in [0, M]$ and $\int_0^M |h(x)|dx < \infty$. Since $M$ is arbitrary, the claim follows. □

Now, fix
$$T \in (0, \infty], \ \gamma \in (0, 1).$$
Note that we allow $T = \infty$. As mentioned before, we will drop the super-indices $T$ in $\mathscr{E}^T$ and $I_\alpha^T$. Suppose $f$ is a distribution in $\mathscr{D}(-\infty, T)$ supported on $[0, T)$. We then formally denote $g_\gamma * f$ which should again be understood as in Equation (17) by $\frac{1}{\Gamma(\gamma)} \int_0^t (t-s)^{\gamma-1} f(s) ds, t \in [0, T)$. We say $f$ is locally integrable function if we can find a locally integrable function $\tilde{f}$ such that $\langle f, \varphi \rangle = \int \tilde{f} \varphi dt, \forall \varphi \in C_c^\infty(-\infty, T)$. It is almost trivial to see that

LEMMA 3.2. *Let $\gamma \in (0, 1)$. If $f \in L^1_{\text{loc}}[0, T)$,*

$$(28) \qquad J_\gamma(f)(t) = g_\gamma * (\theta f)(t) = \frac{1}{\Gamma(\gamma)} \int_0^t (t-s)^{\gamma-1} f(s) ds, \quad t \in [0, T),$$

*where the integral on the right is understood in the Lebesgue sense.*

We introduce

$$(29) \qquad X := \left\{ \varphi \in L^1_{\text{loc}}[0, T) : \exists \varphi_0 \in \mathbb{R}, \lim_{t \to 0^+} \frac{1}{t} \int_0^t |\varphi - \varphi_0| dt = 0 \right\}.$$

The space $X$ is indeed the space of locally integrable functions with $t = 0$ being a Lebesgue point from the right (with the value at $t = 0$ defined as $\varphi_0$). Clearly, $X$ is a vector space and $C[0, T) \subset X \subset L^1_{\text{loc}}[0, T)$. The value $\varphi_0$ is unique for every $\varphi \in X$. We denote

$$(30) \qquad \varphi(0+) := \varphi_0.$$

The notation $\varphi(0+)$ is understood as the right limit in a weak sense and we do not necessarily have a Lebesgue measure zero set $\mathscr{N}$ such that $\lim_{t \to 0^+, t \notin \mathscr{N}} |\varphi(t) - \varphi_0| = 0$.

For convenience, we also introduce the following set for $0 < \gamma < 1$:

$$(31) \qquad Y_\gamma := \left\{ f \in L^1_{\text{loc}}[0, T) : \lim_{T \to 0^+} \frac{1}{T} \int_0^T \left| \int_0^t (t-s)^{\gamma-1} f(s) ds \right| dt = 0 \right\},$$

and also

$$(32) \qquad X_\gamma := \left\{ \varphi : \exists C \in \mathbb{R}, f \in Y_\gamma, \text{s.t. } \varphi = C + J_\gamma(f) \right\}.$$

Recall that for $t \in [0, T)$, we have
$$J_\gamma(f)(t) = g_\gamma * (\theta f)(t) = \frac{1}{\Gamma(\gamma)} \int_0^t (t-s)^{\gamma-1} f(s) ds,$$
where the integral is in the Lebesgue sense by Lemma 3.2.

By the definition of $Y_\gamma$, it is almost trivial to conclude that:

LEMMA 3.3. *$Y_\gamma$ and $X_\gamma$ are subspaces of $L^1_{\text{loc}}[0, T)$. If $f \in Y_\gamma$, then $J_\gamma f(0+) = 0$ and $X_\gamma \subset X$.*

REMARK 3. *If $f \geq 0, a.e.$, $\lim_{t \to 0+} \frac{1}{t} \int_0^t |\int_0^\tau (\tau - s)^{\gamma-1} f(s) ds| d\tau = 0$ is equivalent to $\lim_{t \to 0+} \frac{1}{t} \int_0^t (t-s)^\gamma f(s) ds = 0$. Hence, $L_{loc}^{1/\gamma}[0,T) \subset Y_\gamma$. (e.g. $t^{-\gamma} \notin Y_\gamma$ while $t^{-\gamma+\delta} \in Y_\gamma, \forall \delta > 0$.)*

Now, we introduce our definition of Caputo derivatives based on the modified Riemann–Liouville operators by removing the singularity due to the jump at $t = 0$:

DEFINITION 3.4. *Let $0 < \gamma < 1$ and $\varphi \in L_{loc}^1[0,T)$. For any $\varphi_0 \in \mathbb{R}$, we define the generalized Caputo derivative associated with $\varphi_0$ to be:*

$$D_c^\gamma : \varphi \mapsto D_c^\gamma \varphi = J_{-\gamma} \varphi - \varphi_0 g_{1-\gamma} = J_{-\gamma}(\varphi - \varphi_0) \in \mathscr{E}. \tag{33}$$

*If $\varphi \in X$, we impose $\varphi_0 = \varphi(0+)$ unless explicitly mentioned and in this case, we call $D_c^\gamma$ the Caputo derivative of order $\gamma$.*

Note that if $\varphi$ does not have regularities, $D_c^\gamma \varphi$ is generally a distribution in $\mathscr{E}$. If $\varphi_0 = 0$, it is clear that the generalized Caputo derivative is the Riemann–Liouville derivative. We impose $\varphi_0 = \varphi(0+)$ whenever $\varphi \in X$ to remove the singularity introduced by the jump at $t = 0$, and to be consistent with the traditional Caputo derivative (3) (see Proposition 3.6). This convention is convenient for studying the initial value problems in Section 4. We also point out that this generalized definition can be applied for Caputo derivatives from any time point $t = a$ (If $a = -\infty$, we just remove the limit at $-\infty$). One only uses operators $I_\alpha(\theta(\cdot - a)\varphi)$ to replace $J$.

Recall $J_{-\gamma}\varphi = g_{-\gamma} * (\theta\varphi)$ and in the case $T < \infty$, it is understood as in Equations (17). Note that we have used explicitly the convolution operator $J_{-\gamma}$ in the definition. The convolution structure here enables us to establish the fundamental theorem (Theorem 3.7) below using deconvolution so that we can rewrite fractional differential equations using integral equations with completely monotone kernels.

LEMMA 3.5. *By the definition, we have the following claims:*
1. *For any constant $C$, $D_c^\gamma C = 0$.*
2. *$D_c^\gamma : X \to \mathscr{E}$ is a linear operator.*
3. *$\forall \varphi \in X$, $0 < \gamma_1 < 1$ and $\gamma_2 > \gamma_1 - 1$, we have*

$$J_{\gamma_2} D_c^{\gamma_1} \varphi = \begin{cases} D_c^{\gamma_1-\gamma_2} \varphi, & \gamma_2 < \gamma_1, \\ J_{\gamma_2-\gamma_1}(\varphi - \varphi(0+)), & \gamma_2 \geq \gamma_1. \end{cases}$$

4. *Suppose $0 < \gamma_1 < 1$. If $f \in Y_{\gamma_1}$, $D_c^{\gamma_2} J_{\gamma_1} f = J_{\gamma_1-\gamma_2} f$ for $0 < \gamma_2 < 1$.*
5. *If $D_c^{\gamma_1} \varphi \in X$, then for $0 < \gamma_2 < 1$, $0 < \gamma_1 + \gamma_2 < 1$,*

$$D_c^{\gamma_2} D_c^{\gamma_1} \varphi = D_c^{\gamma_1+\gamma_2} \varphi - D_c^{\gamma_1}\varphi(0+) g_{1-\gamma_2}.$$

6. *$J_{\gamma-1} D_c^\gamma \varphi = J_{-1}\varphi - \varphi(0+)\delta(t)$. If we define this to be $D_c^1$, then for $\varphi \in C^1[0,T)$, $D_c^1 \varphi = \varphi'$.*

*Proof.* The first follows from $g_{-\gamma} * \theta(t) = g_{-\gamma} * g_1 = g_{1-\gamma}$. The second is obvious. The third claim follows from $J_{\gamma_2} D_c^{\gamma_1}\varphi = J_{\gamma_2}(J_{-\gamma_1}\varphi - \varphi(0+)g_{1-\gamma_1}) = J_{\gamma_2-\gamma_1}\varphi - \varphi(0+)g_{1-\gamma_1+\gamma_2}$, which holds by the group property. For the fourth, we just note that $J_{\gamma_1} f(0+) = 0$ and use the group property for $J_\alpha$. The fifth statement follows easily from the third statement. The last claim follows from $J_{\gamma-1}(J_{-\gamma}\varphi - \varphi(0+)g_{1-\gamma}) = J_{-1}\varphi - \varphi(0+)g_0$ and Equation (22). □

Now, we verify that our definition agrees with (3) if $\varphi$ has some regularity. For $\varphi \in L_{loc}^1(0,T)$, we say $\tilde{\varphi}$ is a version of $\varphi$ if they are different only on a Lebesgue measure zero set.

PROPOSITION 3.6. *For $\varphi \in X$, if the distributional derivative of $\varphi$ on $(0, T)$ is integrable on $(0, T_1)$ for any $T_1 \in (0, T)$, then there is a version of $\varphi$ (still denoted as $\varphi$) that is absolutely continuous on $(0, T_1)$ for any $T_1 < T$, so that*

$$(34) \qquad D_c^\gamma \varphi(t) = \frac{1}{\Gamma(1-\gamma)} \int_0^t \frac{\varphi'(s)}{(t-s)^\gamma} ds, \ t < T,$$

*where the convolution integral can be understood in the Lebesgue sense. Further, $D_c^\gamma \varphi \in L_{\text{loc}}^1[0, T)$.*

*Proof.* The existence of the version of $\varphi$ that is absolutely continuous on $(0, T_1)$, $T_1 \in (0, T)$ is a classical result (see [24]). Consequently, the usual derivative $\varphi'$ is the distributional derivative of $\varphi$ on $(0, T)$.

We first consider $T = \infty$.

Define $\varphi^\epsilon = \varphi * \eta_\epsilon$ where $\eta_\epsilon = \frac{1}{\epsilon} \eta(\frac{t}{\epsilon})$ and $0 \le \eta \le 1$ satisfies: (i). $\eta \in C_c^\infty(\mathbb{R})$ with $\text{supp}(\eta) \subset (-M, 0)$ for some $M > 0$. (ii). $\int \eta dt = 1$. $\varphi^\epsilon$ is clearly smooth. Then, we have in $\mathscr{D}'(\mathbb{R})$

$$D(\theta(t)\varphi^\epsilon)(t) = \delta(t)\varphi^\epsilon(0) + \theta(t)D(\varphi^\epsilon),$$

which can be verified easily. Now, take $\epsilon \to 0$. Let $\varphi_0 = \varphi(0+)$. Then, $|\varphi^\epsilon(0) - \varphi_0| = |\int_0^M \eta(-x)\varphi(\epsilon x)dx - \varphi_0| \le \sup |\eta| \int_0^M |\varphi(\epsilon x) - \varphi_0|dx = \sup|\eta| \frac{M}{\epsilon M} \int_0^{M\epsilon} |\varphi(y) - \varphi_0|dy \to 0$. Hence, $\varphi^\epsilon(0) \to \varphi_0$.

Since $\theta(t)\varphi^\epsilon \to \theta(t)\varphi$ in $L_{\text{loc}}^1(\mathbb{R})$, $D(\theta(t)\varphi^\epsilon) \to D(\theta(t)\varphi)$ in $\mathscr{D}'(\mathbb{R})$. For $\theta(t)D(\varphi^\epsilon)$, by the convolution property we have $D(\varphi^\epsilon) = (D\varphi)^\epsilon$ in $\mathscr{D}'(\mathbb{R})$. Since $\varphi'$ is integrable on $(0, T_1)$ for $T_1 \in (0, \infty)$, we can define its values to be zero on $(-\infty, 0]$ and then it becomes a distribution in $\mathscr{D}'(\mathbb{R})$. We still denote it as $\varphi'$. If $\text{supp}(\eta) \subset (-M, 0)$, then $\theta(t)(D\varphi)^\epsilon \to \theta(t)\varphi'$ in $\mathscr{D}'(\mathbb{R})$. (Take a special notice that if $\varphi$ is a causal function, $D\varphi$ as a distribution in $\mathscr{D}'(\mathbb{R})$ generally has an atom $\delta(t)$, but $\theta(t)\varphi'$ does not include this singularity.) This then verifies the distribution identity:

$$D(\theta(t)\varphi(t)) = \delta(t)\varphi(0+) + \theta(t)\varphi'(t).$$

By the definition of $J_{-\gamma}$ (Definition 2.14) and applying Lemma 2.2,

$$J_{-\gamma}\varphi(t) = \frac{1}{\Gamma(1-\gamma)} D\left(\theta(t)t^{-\gamma}\right) * (\theta(t)\varphi) = \frac{1}{\Gamma(1-\gamma)} \left(\theta(t)t^{-\gamma}\right) * D(\theta(t)\varphi)$$

$$= \frac{1}{\Gamma(1-\gamma)} \left(\theta(t)t^{-\gamma}\right) * (\delta(t)\varphi(0+) + \theta(t)\varphi').$$

The first term gives $\frac{\varphi(0+)}{\Gamma(1-\gamma)} \theta(t) t^{-\gamma}$.

Consider now $(\theta(t)t^{-\gamma}) * (\theta(t)\varphi')$. It is clear that

$$(\theta(t)t^{-\gamma}) * (\theta(t)\varphi') = \lim_{M \to \infty} (\theta(t)t^{-\gamma}\chi(t \le M)) * (\theta(t)\varphi'\chi(t \le M)), \ \text{in} \ \mathscr{D}',$$

where $\chi(E)$ is the indicator function of the set $E$. With the truncation, the two functions on the right-hand side are in $L^1(\mathbb{R})$. The convolution $h_M = (\theta(t)t^{-\gamma}\chi(t \le M)) * (\theta(t)\varphi'\chi(t \le M)) \in L^1(\mathbb{R})$ and

$$h_M(t) = \int_0^t \frac{\varphi'(s)}{(t-s)^\gamma} ds, \ 0 < t < M,$$

where the integral is in the Lebesgue sense. Hence, as $M \to \infty$, we find that $(\theta(t)t^{-\gamma}) * (\theta(t)\varphi')$ is a measurable function and for almost every $t$, Equation (34) holds and the integral is a Lebesgue integral. By Lemma 3.1, $D_c^\gamma \varphi \in L_{\text{loc}}^1[0, \infty)$.

Then, by Definition 3.4, we obtain Equation (34). This then finishes the proof for $T = \infty$.

For $T < \infty$, consider $K_n^T \varphi = \chi_n \varphi$ is absolutely continuous on $(0, \infty)$. By the result just proved, we have

$$D_c^\gamma (K_n^T \varphi)(t) = \frac{1}{\Gamma(1-\gamma)} \int_0^t \frac{\chi_n' \varphi + \chi_n \varphi'}{(t-s)^\gamma} ds.$$

It follows that for $t < T - \frac{1}{n}$,

$$D_c^\gamma (K_n^T \varphi)(t) = \frac{1}{\Gamma(1-\gamma)} \int_0^t \frac{\varphi'}{(t-s)^\gamma} ds.$$

Hence, in $\varphi \in \mathscr{D}'(-\infty, T)$, $D_c^\gamma \varphi = \lim_{n \to \infty} R^T D_c^\gamma (K_n^T \varphi)$ is given by the formula listed in the statement. □

Regarding the Caputo derivatives, we introduce several results that may be applied for fractional ODEs and fractional PDEs. The first is the fundamental theorem of fractional calculus:

THEOREM 3.7. *Suppose $\varphi \in L^1_{\text{loc}}[0, T)$ and denote $f = D_c^\gamma \varphi$ $(0 < \gamma < 1)$ to be the generalized Caputo derivative associated with $\varphi_0$, supported in $[0, T)$. Then, for Lebesgue a.e. $t \in (0, T)$,*

(35) $$\varphi(t) = \varphi_0 + J_\gamma(f)(t).$$

*If $f \in L^1_{\text{loc}}[0, T)$, we have for Lebesgue a.e. $t \in (0, T)$ that*

(36) $$\varphi(t) = \varphi_0 + \frac{1}{\Gamma(\gamma)} \int_0^t (t-s)^{\gamma-1} f(s) ds,$$

*where the integral is understood in Lebesgue integral sense.*

*Proof.* By the definition (Equation (33)), $f(t) + \varphi_0 \frac{\theta(t)}{\Gamma(1-\gamma)} t^{-\gamma} = J_{-\gamma}(\varphi(t))$. Then, the convolution group property yields:

$$\theta(t)\varphi(t) = J_\gamma \left( f + \varphi_0 \frac{\theta(t)}{\Gamma(1-\gamma)} t^{-\gamma} \right) = J_\gamma f + \frac{\theta(t)\varphi_0}{\Gamma(\gamma)\Gamma(1-\gamma)} \int_0^t (t-s)^{\gamma-1} s^{-\gamma} ds.$$

Since $\int_0^t (t-s)^{\gamma-1} s^{-\gamma} ds = B(\gamma, 1-\gamma) = \Gamma(\gamma)\Gamma(1-\gamma)$, the second term is just $\theta(t)\varphi_0$.

Since $J_\gamma f = \theta(t)\varphi(t) - \theta(t)\varphi_0 \in L^1_{\text{loc}}[0, T)$, the equality in the distributional sense implies equality in $L^1_{\text{loc}}[0, T)$ and thus a.e.. Further, if $f \in L^1_{\text{loc}}[0, T)$, the second part is obvious by Lemma 3.2. □

This theorem is fundamental for fractional differential equations because this allows us to transform the fractional differential equations to integral equations with completely monotone kernels (which are nonnegative). Then, we are able to establish the comparison principle (Theorem 4.10), good for *a priori* estimate of fractional PDEs.

Using Theorem 3.7, we conclude that

COROLLARY 3.8. *Suppose $\varphi(t) \in L^1_{\text{loc}}[0, T)$, $\varphi \geq 0$ and $\varphi_0 = 0$. If the generalized Caputo derivative $D_c^\gamma \varphi$ $(0 < \gamma < 1)$ associated with $\varphi_0$ is locally integrable, and $D_c^\gamma \varphi \leq 0$, then $\varphi(t) = 0$. (The local integrability assumption can be dropped if we understand the inequality in the distribution sense as in Section 4.)*

Now, we consider functions whose Caputo derivatives are $L^1_{\text{loc}}[0,T)$. Recall the definitions of $X$ in (29), $X_\gamma$ in (32) and $Y_\gamma$ in (31).

PROPOSITION 3.9. *Let $f \in L^1_{\text{loc}}[0,T)$. Then, $D^\gamma_c \varphi = f$ has solutions $\varphi \in X$ if and only if $f \in Y_\gamma$. If $f \in Y_\gamma$, the solutions are in $X_\gamma$ and they can be written as*

$$\varphi(t) = C + \frac{1}{\Gamma(\gamma)} \int_0^t (t-s)^{\gamma-1} f(s) ds, \quad \forall C \in \mathbb{R}.$$

*Further, $\forall \varphi \in X_\gamma$, $D^\gamma_c \varphi \in Y_\gamma$.*

*Proof.* Suppose that $D^\gamma_c \varphi = f$ has a solution $\varphi \in X$. Since $f \in L^1_{\text{loc}}[0,T)$, by Theorem 3.7, we have

$$\varphi(t) = \varphi(0+) + \frac{1}{\Gamma(\gamma)} \int_0^t (t-s)^{\gamma-1} f_1(s) ds,$$

and the integral is in the Lebesgue sense. Since $\varphi \in X$, we have $\lim_{T \to 0} \frac{1}{T} \int_0^T |\varphi(s) - \varphi(0+)| ds = 0$. It follows that

$$\frac{1}{T} \int_0^T \left| \int_0^t (t-s)^{\gamma-1} f_1(s) ds \right| dt \to 0, \quad T \to 0.$$

Hence, $f_1 \in Y_\gamma$. This implies that there are no solutions for $D^\gamma_c \varphi = f$ if $f \in L^1_{\text{loc}}[0,T) \setminus Y_\gamma$. (For example, $D^\gamma_c \varphi = t^{-\gamma}$ has no solutions in $X$.)

For the other direction, now assume $f \in Y_\gamma$. We first note $D^\gamma_c \varphi = 0$ implies that $\varphi$ is a constant by Theorem 3.7. One then can check that $J_\gamma f$, which is in $X_\gamma$ by definition, is a solution to the equation $D^\gamma_c \varphi = f$. Hence any solution can be written as $J_\gamma f + C$. The other direction and the second claim are shown.

We now show the last claim. If $\varphi \in X_\gamma$, by definition (Equation (32)), $\exists f \in Y_\gamma, C \in \mathbb{R}$ such that

$$\varphi(t) = C + \frac{1}{\Gamma(\gamma)} \int_0^t (t-s)^{\gamma-1} f(s) ds.$$

Since $f \in Y_\gamma$, $J_\gamma f(0+) = 0$ by Lemma 3.3. This means $C = \varphi(0+)$. Now, apply $J_{-\gamma}$ on both sides. Note that $J_{-\gamma} C = C g_{-\gamma} * g_1 = C g_{1-\gamma} = C \frac{\theta(t)}{\Gamma(1-\gamma)} t^{-\gamma}$. By the group property, we find that

$$D^\gamma_c \varphi = J_{-\gamma} J_\gamma f = f. \qquad \square$$

Motivated by the discussion in Section 2.4, we have another claim about the equation $D^\gamma_c \varphi = f$ where the solutions are in $C^0[0,T)$:

PROPOSITION 3.10. *Suppose $f \in L^1[0,T) \cap \tilde{H}^s(0,T)$. If $s$ satisfies (i). $s \geq 0$ when $\gamma > 1/2$ or (ii). $s > \frac{1}{2} - \gamma$ when $\gamma \leq 1/2$, then $\exists \varphi \in C^0[0,T)$ such that $D^\gamma_c \varphi = f$. If $T = \infty$, we can also ask for $f \in L^1_{\text{loc}}[0,\infty) \cap \tilde{H}^s_{\text{loc}}(0,\infty)$.*

Let us focus on the mollifying effect on the Caputo derivatives. Let $\eta \in C^\infty_c(\mathbb{R})$, $0 \leq \eta \leq 1$ and $\int \eta \, dt = 1$. We define $\eta_\epsilon = \frac{1}{\epsilon} \eta(\frac{t}{\epsilon})$. Consider $\varphi \in L^1[0,\infty)$, it is well known that

(37) $$\varphi^\epsilon = \varphi * \eta_\epsilon \in C^\infty(\mathbb{R})$$

and that $\text{supp}(\varphi^\epsilon) \subset \text{supp}(\varphi) + \text{supp}(\eta_\epsilon)$ where $A + B = \{x + y : x \in A, y \in B\}$.

PROPOSITION 3.11. *Suppose $T = \infty$. Assume $\mathrm{supp}(\eta) \subset (-\infty, 0)$.*

*(i) $\forall \varphi \in X$, $D_c^\gamma(\varphi^\epsilon) \to D_c^\gamma \varphi$ in $\mathscr{D}'(\mathbb{R})$ as $\epsilon \to 0^+$. Also,*

$$\text{(38)} \quad D_c^\gamma \varphi^\epsilon(t) = \frac{1}{\Gamma(1-\gamma)} \int_0^t \frac{(\varphi^\epsilon)'(s)}{(t-s)^\gamma} ds$$
$$= \frac{1}{\Gamma(1-\gamma)} \left( \frac{\varphi^\epsilon(t) - \varphi^\epsilon(0)}{t^\gamma} + \gamma \int_0^t \frac{\varphi^\epsilon(t) - \varphi^\epsilon(s)}{(t-s)^{\gamma+1}} ds \right).$$

*(ii) If $E(\cdot)$ is a $C^1$ convex function, then*

$$\text{(39)} \quad D_c^\gamma E(\varphi^\epsilon) \leq E'(\varphi^\epsilon) D_c^\gamma \varphi^\epsilon.$$

*If there exists a sequence $\epsilon_k$, such that $D_c^\gamma \varphi^{\epsilon_k}$ converges in $L_{\mathrm{loc}}^1[0, \infty)$. Then, the limit is $D_c^\gamma \varphi$ and $D_c^\gamma \varphi \in L_{\mathrm{loc}}^1[0, \infty)$. Moreover, in the distributional sense,*

$$\text{(40)} \quad D_c^\gamma E(\varphi) \leq E'(\varphi) D_c^\gamma \varphi.$$

*(iii) If $\varphi \in C[0, T_1] \cap C^1(0, T_1]$ for some $T_1 > 0$, we have for all $t \in (0, T_1]$*

$$\text{(41)} \quad D_c^\gamma E(\varphi)(t) \leq E'(\varphi(t)) D_c^\gamma \varphi(t).$$

*Proof.* (i). For any $\varphi \in X$, it is clear that $\theta(t)\varphi^\epsilon \to \theta(t)\varphi$ in $L_{\mathrm{loc}}^1[0, \infty)$ and hence in the distributional sense. Using the definition of $J_\alpha$ and the definition of convolution on $\mathscr{E}$, one can readily check $J_{-\gamma} \varphi^\epsilon \to J_{-\gamma} \varphi$ in $\mathscr{D}'(\mathbb{R})$.

For $\varphi^\epsilon(0) \to \varphi(0+)$, we need $\mathrm{supp}(\eta) \subset (-\infty, 0)$. There then exists $M > 0$ such that $\eta(t) = 0$ if $t < -M$. Let $\varphi_0 = \varphi(0+)$. Then, $|\varphi^\epsilon(0) - \varphi_0| = |\int_0^M \eta(-x)\varphi(\epsilon x)dx - \varphi_0| \leq \sup |\eta| \int_0^M |\varphi(\epsilon x) - \varphi_0| dx = \sup |\eta| \frac{M}{\epsilon M} \int_0^{M\epsilon} |\varphi(y) - \varphi_0| dy \to 0$. Hence, $\varphi^\epsilon(0) \to \varphi_0$. This then shows that the first claim is true.

For the alternative expressions of $D_c^\gamma \varphi^\epsilon$, we have used Proposition 3.6 and integration by parts. These are valid since $\varphi^\epsilon \in C^\infty$.

(ii). Multiplying $E'(\varphi^\epsilon(t))$ on $\frac{\varphi^\epsilon(t) - \varphi^\epsilon(0)}{t^\gamma} + \gamma \int_0^t \frac{\varphi^\epsilon(t) - \varphi^\epsilon(s)}{(t-s)^{\gamma+1}} ds$ and using the inequality

$$E'(a)(a - b) \geq E(a) - E(b)$$

since $E(\cdot)$ is convex, we find

$$\Gamma(1-\gamma) E'(\varphi^\epsilon(t)) D_c^\gamma \varphi^\epsilon(t) \geq \frac{E(\varphi^\epsilon(t)) - E(\varphi^\epsilon(0))}{t^\gamma} + \gamma \int_0^t \frac{E(\varphi^\epsilon(t)) - E(\varphi^\epsilon(s))}{(t-s)^{\gamma+1}} ds.$$

Since $E$ is $C^1$ and $\varphi^\epsilon$ is smooth, the second integral converges for almost every $t$. Further, $E(\varphi^\epsilon)' = E'(\varphi^\epsilon)(\varphi^\epsilon)'$ a.e.. Then, the right-hand side must be $\int_0^t \frac{E(\varphi^\epsilon)'}{(t-s)^\gamma} ds$. By Proposition 3.6, we end the second claim.

The last claim is trivial by sending $\epsilon_k \to 0$ in (39) and using the $L_{\mathrm{loc}}^1$ convergence of Caputo derivative (the inequality will be preserved under the limit even in the distributional sense).

(iii). We extend $\varphi$ such that it is in $C[0, \infty) \cap C^1(0, \infty)$. Then, $\varphi$ is absolutely continuous on any closed interval. It is not hard to see $D_c^\gamma \varphi^{\epsilon_k}$ converges in $L_{\mathrm{loc}}^1[0, \infty)$ to $D_c^\gamma \varphi$ using Proposition 3.6 and the regularity of $\varphi$. Hence, by (ii), we have

$$D_c^\gamma E(\varphi) \leq E'(\varphi) D_c^\gamma \varphi,$$

in the distributional sense. However, by Proposition 3.6, all the functions here are in $C(0, \infty)$. Noting that the extension does not change the derivatives on $(0, T_1]$, the desired claim follows. □

REMARK 4. *It is interesting to note that we have to choose $\eta$ such that $\mathrm{supp}(\eta) \subset (-\infty, 0)$. Using other mollifiers, we may not get the correct limit. This reflects that Caputo derivatives only model the dynamics of memory from $t = 0^+$ and the singularities at $t = 0$ for Riemann–Liouville derivatives are removed totally. It is exactly this nature that makes Caputo derivatives to have many similarities with the ordinary derivative and suitable for initial value problems.*

Now, as consequences of Proposition 3.11, we verify that our definition of Caputo derivative can recover (4) as used in [1] by Allen, Caffarelli and Vasseur. We use $C^\alpha[0,T]$ to mean the set of functions that are $\alpha$-Hölder continuous (see [9, Chap. 5]) on $[0,T]$. We have the following:

COROLLARY 3.12. *Suppose $T \in (0,\infty]$. If there exists $\delta > 0$ such that $\varphi \in C^{\gamma+\delta}[0,T)$, then $D_c^\gamma \varphi \in C[0,T)$ and $\forall t \in [0,T)$,*

$$D_c^\gamma \varphi(t) = \frac{1}{\Gamma(1-\gamma)} \left( \frac{\varphi(t) - \varphi(0)}{t^\gamma} + \gamma \int_0^t \frac{\varphi(t) - \varphi(s)}{(t-s)^{\gamma+1}} ds \right).$$

*Proof.* We show the claim for $T < \infty$. The proof for $T = \infty$ is similar but easier. As in Equation (14), we define $\varphi_n = K_n^T \varphi = \chi_n \varphi$. $\varphi_n \in L^1[0,\infty)$ and $\varphi_n = \varphi$ in $[0, T-1/n]$ and $\varphi_n$ is also $\gamma + \delta$-Hölder continuous. Pick the mollifier $\eta$ such that $\mathrm{supp}\, \eta \subset (-\infty, 0)$, and $\varphi_n^\epsilon = \eta_\epsilon * \varphi_n$. Then, $\mathrm{supp}\, \varphi_n^\epsilon \subset (-\infty, T)$. By Proposition 3.11,

$$D_c^\gamma \tilde{\varphi}_n^\epsilon(t) = \frac{1}{\Gamma(1-\gamma)} \left( \frac{\varphi_n^\epsilon(t) - \varphi_n^\epsilon(0)}{t^\gamma} + \gamma \int_0^t \frac{\varphi_n^\epsilon(t) - \varphi_n^\epsilon(s)}{(t-s)^{\gamma+1}} ds \right).$$

By the $\gamma + \delta$-Hölder continuity, $D_c^\gamma \tilde{\varphi}_n^\epsilon$ converges uniformly on $[0, T-1/n]$ to

$$\frac{1}{\Gamma(1-\gamma)} \left( \frac{\varphi_n(t) - \varphi_n(0)}{t^\gamma} + \gamma \int_0^t \frac{\varphi_n(t) - \varphi_n(s)}{(t-s)^{\gamma+1}} ds \right),$$

which must be $D_c^\gamma \varphi_n$ since $\varphi_n^\epsilon \to \varphi_n$ uniformly by Proposition 3.11. The uniform limit must be continuous too. By the definition of $K_n^T$, for any $t < T - 1/n$,

$$D_c^\gamma \varphi_n(t) = \frac{1}{\Gamma(1-\gamma)} \left( \frac{\varphi(t) - \varphi(0)}{t^\gamma} + \gamma \int_0^t \frac{\varphi(t) - \varphi(s)}{(t-s)^{\gamma+1}} ds \right).$$

Therefore, since the limit of $D_c^\gamma \varphi_n$ in $\mathscr{D}'(-\infty, T)$ is $D_c^\gamma \varphi$, $D_c^\gamma \varphi$ must be given as in the statement. □

Lastly, we consider the Laplace transform of the generalized Caputo derivatives in the case $T = \infty$. The generalized Caputo derivative $D_c^\gamma \varphi$ associated with $\varphi_0$ is only a distribution. Recalling that its support is in $[0,\infty)$, we then define the Laplace transform of $D_c^\gamma \varphi$ as

(42) $$\mathcal{L}(D_c^\gamma \varphi) = \lim_{M \to \infty} \langle D_c^\gamma \varphi, \zeta_M e^{-st} \rangle,$$

where $\zeta_M(t) = \zeta_0(t/M)$. $\zeta_0 \in C_c^\infty, 0 \leq \zeta_0 \leq 1$ satisfies: (i) $\mathrm{supp}\, \zeta_0 \subset [-2,2]$ (ii) $\zeta_0 = 1$ for $t \in [-1,1]$. This definition clearly agrees with the usual definition of Laplace transform if the usual Laplace transform of function $\varphi$ exists.

We introduce the following set

(43) $$\mathcal{E}(\mathcal{L}) := \left\{ \varphi \in L^1_{\mathrm{loc}}[0,\infty) : \exists L > 0, s.t. \lim_{A \to \infty} \|e^{-Lt}\varphi\|_{L^\infty[A,\infty)} = 0 \right\}.$$

PROPOSITION 3.13. *If $\varphi \in \mathcal{E}(\mathcal{L})$, then for any given $\varphi_0 \in \mathbb{R}$ and the generalized Caputo derivative $D_c^\gamma \varphi$ associated with $\varphi_0$, $\mathcal{L}(D_c^\gamma \varphi)$ is defined for $\mathrm{Re}(s) > L$ and is given by*

$$\mathcal{L}(D_c^\gamma \varphi) = \mathcal{L}(\varphi) s^\gamma - \varphi_0 s^{\gamma-1}. \tag{44}$$

*Proof.* $\zeta_M e^{-st} \in C_c^\infty$. Then, it follows that

$$\langle D_c^\gamma \varphi, \zeta_M e^{-st} \rangle = \left\langle g_{-\gamma} * (\theta(t)\varphi) - \frac{\varphi_0 \theta(t) t^{-\gamma}}{\Gamma(1-\gamma)}, \zeta_M e^{-st} \right\rangle$$

$$= -\frac{1}{\Gamma(1-\gamma)} \left\langle (\theta(t) t^{-\gamma}) * (\theta(t)\varphi), \zeta_M' e^{-st} - s\zeta_M e^{-st} \right\rangle - \frac{\varphi_0}{\Gamma(1-\gamma)} \int_0^\infty t^{-\gamma} \zeta_M e^{-st} dt.$$

Note that $\mathrm{supp}\, \zeta_M' \cap [0, \infty) \subset [M, 2M]$:

$$\left| \int_M^{2M} \int_0^t (t-\tau)^{-\gamma} \varphi(\tau) d\tau \zeta_M' e^{-st} dt \right|$$

$$\leq \frac{\sup|\zeta_0'|}{M} \int_M^{2M} \int_0^t (t-\tau)^{-\gamma} e^{-L\tau} |\varphi(\tau)| d\tau e^{-(\mathrm{Re}(s)-L)t} dt.$$

By the assumption, there exists $T_0$ such that $|e^{-L\tau} \varphi(\tau)| < 1, a.e.$ if $\tau > T_0$. Hence, if $M > 2T_0$, the inner integral is controlled by $T_0^{-\gamma} \int_0^{T_0} |\varphi(\tau)| d\tau + \int_{T_0}^t (t-\tau)^{-\gamma} d\tau \leq C(1+t^{1-\gamma})$. Since $\lim_{M \to \infty} \int_M^{2M} (1+t^{1-\gamma}) e^{-\epsilon t} dt = 0$ for any $\epsilon > 0$, we find that the term associated with $\zeta_M'$ tends to zero as $M \to \infty$.

For the second term,

$$\left\langle (\theta(t) t^{-\gamma}) * (\theta(t)\varphi), \zeta_M e^{-st} \right\rangle = \int_0^\infty \int_0^t (t-\tau)^{-\gamma} \varphi(\tau) d\tau \zeta_M(t) e^{-st} dt$$

$$= \int_0^\infty \varphi(\tau) e^{-s\tau} \int_\tau^\infty (t-\tau)^{-\gamma} \zeta_M(t) e^{-(t-\tau)s} dt.$$

As $M \to \infty$, one finds that

$$\int_0^\infty t^{-\gamma} \zeta_M(t+\tau) e^{-ts} dt \to \Gamma(1-\gamma) s^{\gamma-1},$$

for every $\tau > 0$. Since $\mathrm{Re}(s) > L$, the dominate convergence theorem implies that the first term goes to $\mathcal{L}(\varphi) s^\gamma$.

Similarly, the last term converges to $-\varphi_0 s^{\gamma-1}$. $\square$

To conclude, there is no group property for Caputo derivatives. However, the Caputo derivatives remove the singularities at $t = 0$ compared with the Riemann–Liouville derivatives and have many properties that are similar to the ordinary derivative so that they are suitable for initial value problems.

**4. Time fractional ordinary differential equations.** In this section, we prove some results about time fractional ODEs using the Caputo derivatives, whose new definition and properties have been discussed in Section 3. Compared with those in [8, 20, 7], the assumptions here are sufficiently weak and conclusions are general.

Consider the initial value problem (IVP) for $v : U \to \mathbb{R}$

$$D_c^\gamma v(t) = f(t, v(t)), \quad v(0) = v_0, \tag{45}$$

where $U \subset \mathbb{R}$ is some interval to be determined and $f$ is a measurable function. The initial value $v_0$ is understood to be the initial value used in the generalized Caputo derivative $D_c^\gamma v$. (If $v \in X$, we impose $v_0 = v(0+)$ as in (30).)

DEFINITION 4.1. *Let $T > 0$. A function $v \in L^1_{\text{loc}}[0, T)$ is called a weak solution to (45) on $[0, T)$ if $f(t, v(t)) \in \mathscr{D}'(-\infty, T)$ and (45) holds in the distributional sense. If (i) $v$ is a weak solution and $v \in X$ with $v(0+) = v_0$; (ii) both $D_c^\gamma v$ and $f(t, v(t))$ are locally integrable so that (45) holds a.e. with respect to Lebesgue measure, we call $v(\cdot)$ a strong solution on $[0, T)$.*

It is clear that a solution on $[0, T)$ is also a solution on $[0, T_1)$ for any $T_1 \in (0, T)$. Making use of Theorem 3.7, we have the following equivalence claim:

PROPOSITION 4.2. *Suppose $f \in L^\infty_{\text{loc}}([0, \infty) \times \mathbb{R}; \mathbb{R})$. Fix $T > 0$. Then, $v(t) \in L^1_{\text{loc}}[0, T)$ with initial value $v_0$ is a strong solution of (45) on $(0, T)$ if and only if $v \in X$ and it solves the following integral equation*

$$(46) \qquad v(t) = v_0 + \frac{1}{\Gamma(\gamma)} \int_0^t (t-s)^{\gamma-1} f(s, v(s)) ds, \ \forall t \in (0, T).$$

The '$\Rightarrow$' direction follows from Theorem 3.7. For the other direction, if the integral equation holds, then $\int_0^t (t-s)^{\gamma-1} |f(s, v(s))| \, ds < \infty$ a.e. by the definition of Lebesgue integral. Since $(t-s)^{\gamma-1} \geq t^{\gamma-1}$ on $[0, t]$, we therefore know $f(s, v(s))$ is integrable on $[0, t]$ for a.e. $t \in [0, T)$, which further implies that $f(s, v(s))$ is indeed integrable on $[0, T - \delta]$ for any $\delta \in (0, T)$. The group property of $J_\alpha$ ensures that the equation holds in the distributional sense. Hence, $v \in X$ is a strong solution to (45).

We start with a simple linear fractional ODE:

PROPOSITION 4.3. *Let $0 < \gamma < 1$, $\lambda \neq 0$, and suppose $b(t)$ is continuous such that there exists $L > 0$, $\limsup_{t \to \infty} e^{-Lt} |b(t)| = 0$. Then, there is a unique strong solution of the equation*

$$D_c^\gamma v = \lambda v + b, \quad v(0) = v_0$$

*in $\mathcal{E}(\mathcal{L})$ (see Equation (43)) and is given by*

$$(47) \qquad v(t) = v_0 e_{\gamma, \lambda}(t) + \frac{1}{\lambda} \int_0^t b(t-s) e'_{\gamma, \lambda}(s) ds,$$

*where $e_{\gamma, \lambda}(t) = E_\gamma(\lambda t^\gamma)$ and*

$$(48) \qquad E_\gamma(z) = \sum_{n=0}^\infty \frac{z^n}{\Gamma(\gamma n + 1)}$$

*is the Mittag–Leffler function [16].*

*Proof.* For $\gamma \in (0, 1)$, we have by Proposition 3.13

$$\mathcal{L}(D_c^\gamma v) = s^\gamma V(s) - v_0 s^{\gamma - 1}$$

and $V(s) = \mathcal{L}(v)$.

By the assumption on $b(t)$, the Laplace transform $B(s) = \mathcal{L}(b)$ exists. Hence, any continuous solution of the initial value problem that is in $\mathcal{E}(\mathcal{L})$ must satisfy

$$V(s) = v_0 \frac{s^{\gamma - 1}}{s^\gamma - \lambda} + \frac{B(s)}{s^\gamma - \lambda}.$$

Using the equality ([21, Appendix])

$$\int_0^\infty e^{-st} E_\gamma(s^\gamma z t^\gamma) dt = \frac{1}{s(1-z)}, \tag{49}$$

and denoting $e_{\gamma,\lambda}(t) = E_\gamma(\lambda t^\gamma)$, we have that

$$\mathcal{L}(e_{\gamma,\lambda}) = \frac{s^{\gamma-1}}{s^\gamma - \lambda}, \quad \mathcal{L}(e'_{\gamma,\lambda}) = \frac{\lambda}{s^\gamma - \lambda}.$$

Taking the inverse Laplace transform of $V(\cdot)$, we get (47). Note that though $e'_{\gamma,\lambda}$ blows up at $t = 0$, it is integrable near $t = 0$ and the convolution is well defined. Further, by the asymptotic behavior of $b$ and the decaying rate of $e'_{\gamma,\lambda}$, the solution is again in $\mathcal{E}(\mathcal{L})$. The existence part is proved.

Since the Laplace transform of functions that are in $\mathcal{E}(\mathcal{L})$ is unique, the uniqueness part is proved. □

REMARK 5. *For the existence of solutions in $X$ (where $T = \infty$), the condition $\limsup_{t\to\infty} e^{-Lt}|b(t)| = 0$ can be removed, since for any $t > 0$, we can redefine $b$ beyond $t$ so that $\limsup_{t\to\infty} e^{-Lt}|b(t)| = 0$. The value of $v(t)$ keeps unchanged by the redefinition according to Formula (47). Hence, (47) gives a solution for any continuous function $b$ in $X$. The uniqueness in $X$ instead of in $\mathcal{E}(\mathcal{L})$ will be established in Theorem 4.4 below.*

We now consider a general fractional ODE when $f(t, v)$ is a continuous function.

THEOREM 4.4. *Let $0 < \gamma < 1$ and $v_0 \in \mathbb{R}$. Consider IVP (45). If there exist $T > 0, A > 0$ such that $f$ is defined and continuous on $D = [0,T] \times [v_0 - A, v_0 + A]$ such that there exists $L > 0$,*

$$\sup_{0 \le t \le T} |f(t, v_1) - f(t, v_2)| \le L|v_1 - v_2|, \quad \forall v_1, v_2 \in [v_0 - A, v_0 + A].$$

*Then, the IVP has a unique strong solution on $[0, T_1)$, and $T_1$ is given by*

$$T_1 = \min\left\{T, \; \sup\left\{t \ge 0 : \frac{M}{\Gamma(1+\gamma)} t^\gamma E_\gamma(Lt^\gamma) \le A\right\}\right\} > 0, \tag{50}$$

*where $M = \sup_{0 \le t \le T} |f(t, v_0)|$ and*

$$E_\gamma(z) = \sum_{n=0}^\infty \frac{z^n}{\Gamma(n\gamma + 1)}$$

*is the Mittag–Leffler function.*

*Moreover, $v(\cdot) \in C[0, T_1]$. Further, the solution is continuous with respect to the initial value. Indeed, fix $t \in (0, T_1)$ and $\epsilon_0 < T_1 - t$. Then, $\forall \epsilon \in (0, \epsilon_0), \exists \delta_0 > 0$ such that for any $|\delta| \le \delta_0$, the solution of the fractional ODE with initial value $v_0 + \delta$, $v^\delta(\cdot)$, exists on $(0, T_1 - \epsilon_0)$ and*

$$|v^\delta(t) - v(t)| < \epsilon. \tag{51}$$

*Proof.* The proof is just like the proof of the existence and uniqueness theorem for ODEs using Picard iteration.

Consider the sequence constructed by

$$v^n(t) = v_0 + g_\gamma * (\theta(\cdot)f(\cdot, v^{n-1}(\cdot)))(t), \ n \geq 1, \tag{52}$$
$$v^0(t) = v_0.$$

where $t \in [0, T_1]$. Recall that the convolution in principle is understood as in Equation (17) and it can be understood as the Lebesgue integral $\frac{1}{\Gamma(\gamma)} \int_0^t (t-s)^{\gamma-1} f(s, v^{n-1}(s)) ds$ by Lemma 3.2.

Consider $E^n = |v^n - v^{n-1}|$. We then find for $t \in [0, T_1]$,

$$E^1(t) = |g_\gamma * (\theta(t)f(t, v^0))| \leq \frac{M}{\Gamma(1+\gamma)} T_1^\gamma =: M_{T_1}.$$

One can verify that $M_{T_1} \leq A$ by the definition of $T_1$.

Now, we assume that $\sum_{m=1}^{n-1} E^m \leq A$ so that $|v^{n-1} - v_0| \leq A$. We will then show that this is true for $E^n$ as well. With this induction assumption, we can find that for $t \in [0, T_1]$ and $m = 2, 3, \ldots, n$, it holds that

$$E^m(t) \leq L g_\gamma * (\theta E^{m-1})(t).$$

Note that $\theta(t) = g_1(t)$ and by the group property, we find on $[0, T_1]$

$$E^m \leq M_{T_1} L^{m-1} g_{1+(m-1)\gamma}, m = 1, 2, \ldots, n.$$

It then follows that

$$\sum_{m=1}^n E^m(t) < M_{T_1} \sum_{m=1}^\infty L^{m-1} g_{1+(m-1)\gamma}(t) = M_{T_1} E_\gamma(L t^\gamma),$$

where $E_\gamma$ is the Mittag–Leffler function. By the definition of $T_1$, $\sum_{m=1}^n E^m(t) \leq A$ for all $t \in [0, T_1]$. Hence, by induction, we have $(t, v^n(t)) \in D$ for all $t \in [0, T_1]$ and $n \geq 0$. It then follows that $\sum_n |v^n - v^{n-1}|$ converges uniformly on $[0, T_1]$. This shows that $v^n \to v$ uniformly on $[0, T_1]$. $v$ is then continuous. Hence, taking $n \to \infty$,

$$v(t) = v_0 + g_\gamma * (\theta(\cdot)f(\cdot, v(\cdot)))(t), \ \ t \in [0, T_1].$$

This means that $v(\cdot)$ is a solution.

We now show the uniqueness. Suppose both $v_1, v_2$ are strong solutions on $[0, T_1)$. Then, $v_1(t)$ and $v_2(t)$ fall into $[v_0 - A, v_0 + A]$ for $t < T_1$ because $f$ is defined on $D$.

Let $w = v_1 - v_2$. Then, by the linearity of $D_c^\gamma$, we have in the distributional sense that

$$D_c^\gamma w(t) = f(t, v_1(t)) - f(t, v_2(t)).$$

Since $f$ is Lipschitz continuous and both $v_1, v_2 \in L^1_{\text{loc}}[0, T_1)$, $f(\cdot, v_1(\cdot)) - f(\cdot, v_2(\cdot)) \in L^1_{\text{loc}}[0, T_1)$. By Theorem 3.7, we have in the Lebesgue sense that for $t \in (0, T_1)$

$$w(t) = \frac{1}{\Gamma(\gamma)} \int_0^t (t-s)^{\gamma-1} (f(s, v_1(s)) - f(s, v_2(s))) ds.$$

For all $t \leq T_1$,

$$|w(t)| \leq \frac{L}{\Gamma(\gamma)} \int_0^t (t-s)^{\gamma-1} |w(s)| ds = L g_\gamma * (\theta |w|)(t).$$

Since $g_\alpha \geq 0$ when $\alpha > 0$ and $\theta g_\alpha = g_\alpha$, we convolve both sides with $g_{n_0\gamma}$ and have

$$\theta(t)w_1(t) \leq L_T g_\gamma * (\theta w_1)(t),$$

where $\theta w_1 = g_{n_0\gamma} * (\theta|w|) \geq 0$. Since $g_{n_0\gamma}(t) = C_1 t^{n_0\gamma - 1}$, if $n_0$ is large enough, $w_1$ is continuous on $[0, T_1]$.

Then, by iteration and the group property, we have

$$\theta(t)w_1(t) \leq L^n g_{n\gamma} * (\theta w_1)(t) \leq \frac{L^n}{\Gamma(n\gamma + 1)} \sup_{0 \leq t \leq T_1} |w_1|(t) \int_0^t (t-s)^{n\gamma} ds.$$

Since $\Gamma(n\gamma + 1)$ grows exponentially, this tends to zero. Hence, $w_1 = 0$ on $[0, T_1]$. Then, convolving both sides with $g_{-n_0\gamma}$ on $w_1 = g_{n_0\gamma} * (\theta(t)|w|) = 0$, we find $|w| = 0$. Hence, $\theta v_1 = \theta v_2$ in $\mathscr{D}'(-\infty, T_1)$ and therefore $v_1 = v_2$ on $[0, T_1)$.

For continuity on the initial value, we make a change of variables $u = v - a$ for $a \in (v_0 - \delta_0, v_0 + \delta_0)$ with suitable $\delta_0 > 0$. Since the Caputo derivative of a constant is zero, the equation is reduced to

$$D_c^\gamma u(t) = f(t, u(t) + a), \quad u(0) = 0.$$

For this question, once again, construct the sequence $u^n$ like in Equation (52). One first show that the $n$th one $u^n$ is continuous on $a \in (v_0 - \delta_0, v_0 + \delta_0)$. Performing similar argument, $u^{n+1} \to u$ uniformly on $[0, T_1 - \epsilon_0]$. Then, the limit $u$ is continuous on $a$. □

COROLLARY 4.5. *Suppose $f$ is defined and continuous on $[0, \infty) \times \mathbb{R}$. If $\forall T > 0$, there exists $L_T$ such that*

$$\sup_{0 \leq t \leq T} |f(t, v_1) - f(t, v_2)| \leq L_T |v_1 - v_2|, \quad \forall v_1, v_2 \in \mathbb{R},$$

*then the unique continuous solution exists on $[0, \infty)$.*

*Proof.* Using the same techniques as in the proof of Theorem 4.4, one can show that the strong solution (not just locally bounded strong solution) is unique. Fix $T > 0$. By the definition of $T_1$, for any given $v_0$, $A$ in Theorem 4.4 can be chosen arbitrarily large so that $T_1 = T$. Hence, the solution exists on $[0, T)$. Since $T$ is arbitrary, the claim follows. □

Now, we establish the following global behavior of the fractional ODE

PROPOSITION 4.6. *Let $-\infty \leq \alpha < \beta \leq \infty$. Assume $f : [0, \infty) \times (\alpha, \beta) \to \mathbb{R}$ is continuous and locally Lipschitz in the second variable. In other words, $\forall A > 0$ and $K \subset (\alpha, \beta)$ compact, there exists $L_{A,K} > 0$ such that*

(53) $$\sup_{0 \leq t \leq A} |f(t, v_1) - f(t, v_2)| \leq L_{A,K} |v_1 - v_2|, \forall v_1, v_2 \in K.$$

*Let $0 < \gamma < 1$ and $v_0 \in (\alpha, \beta)$. Then, the IVP:*

(54) $$D_c^\gamma v(t) = f(t, v(t)), \quad v(0) = v_0,$$

*has a unique continuous solution $v(\cdot)$ on $[0, T_b)$, where*

$$T_b = \sup\Big\{h > 0 : \text{The solution } v \in C[0, h], \ v(t) \in (\alpha, \beta), \ \forall t \in [0, h)\Big\}.$$

*is the largest time of existence satisfying $T_b \in (0, \infty]$. If $T_b < \infty$, then we have $\liminf_{t \to T_b^-} v(t) = \alpha$, or $\limsup_{t \to T_b^-} v(t) = \beta$.*

*Proof.* By Theorem 4.4, strong solution with $[\liminf_{t\to 0^+} v(t), \limsup_{t\to 0^+} v(t)] \subset (\alpha, \beta)$ is unique. This solution exists locally and is continuous. Hence, we have $T_b > 0$. To finish the proof, we only need to show that if $T_b < \infty$ and

$$\alpha < \liminf_{t\to T_b^-} v(t) \leq \limsup_{t\to T_b^-} v(t) < \beta,$$

then the solution can be extended to a larger interval and therefore we have a contradiction. For convenience of notation, let $k_1 = \liminf_{t\to T_b^-} v(t)$ and $k_2 = \limsup_{t\to T_b^-} v(t)$; then $[k_1, k_2] \subset (\alpha, \beta)$ is compact.

Pick $\delta > 0$ so that $[k_1 - \delta, k_2 + \delta] \subset (\alpha, \beta)$. Define another function $\tilde{f}(t, v)$ so that it agrees with $f(t, v)$ on $[0, T_b + \delta] \times [k_1 - \delta, k_2 + \delta]$ and globally Lipschitz. For example, one can choose

$$\tilde{f}(t,v) = \begin{cases} f(t,v) & (t,v) \in [0, T_b+\delta] \times [k_1-\delta, k_2+\delta], \\ f(t, k_2+\delta) & t \leq T_b+\delta, v \geq k_2+\delta, \\ f(T_b+\delta, k_2+\delta) & t \geq T_b+\delta, v \geq k_2+\delta, \\ f(T_b+\delta, v) & t \geq T_b+\delta, v \in [k_1-\delta, k_2+\delta], \\ f(t, k_1-\delta) & t \leq T_b+\delta, v \leq k_1-\delta, \\ f(T_b+\delta, k_1-\delta) & t \geq T_b+\delta, v \leq k_1-\delta. \end{cases}$$

By Corollary 4.5, the unique continuous solution to the time fractional ODE $D_c^\gamma \tilde{v} = \tilde{f}(\cdot, \tilde{v})$ exists on $[0, \infty)$. Since $f$ and $\tilde{f}$ agree on $[0, T_b + \delta] \times [k_1 - \delta, k_2 + \delta]$, $v(\cdot)$ solves $D_c^\gamma \tilde{v} = \tilde{f}(\cdot, \tilde{v})$ as well on $[0, T_b)$. Hence $\tilde{v} = v$ on this interval. It follows that $\exists \delta_1 \in (0, \delta)$ such that $(t, \tilde{v}(t)) \in [0, T_b + \delta] \times [k_1 - \delta, k_2 + \delta]$ for any $t \leq T_b + \delta_1$. Hence, on $[0, T_b + \delta_1]$, $\tilde{v}$ solves $D_c^\gamma v = f(\cdot, v)$ as well, which contradicts with the definition of $T_b$. □

REMARK 6. *Note that we cannot apply the standard continuation technique of ODEs for fractional ODEs, because time fractional ODEs are non-Markovian. In other words, suppose $v_1(t)$ solves the time fractional ODE on $[0, T_1]$ and $v_2(t)$ solves the time fractional ODE with initial value $v_2(0) = v_1(T_1)$ on $[0, \delta]$. Then, the concatenation of $v_1$ and $v_2$ is not a solution to $D_c^\gamma v(t) = f(t, v(t))$ on $[0, T_1 + \delta]$.*

Before further discussion, we introduce the notion of inequalities for distributions.

DEFINITION 4.7. *We say $f \in \mathscr{D}'(-\infty, T)$ is a nonpositive distribution if for any $\varphi \in C_c^\infty(-\infty, T)$ with $\varphi \geq 0$, we have*

(55) $$\langle f, \varphi \rangle \leq 0.$$

*We say $f_1 \leq f_2$ for $f_1, f_2 \in \mathscr{D}'(-\infty, T)$ if $f_1 - f_2$ is nonpositive. We say $f_1 \geq f_2$ if $f_2 - f_1$ is non-positive.*

The following lemma is well known and we omit the proof:

LEMMA 4.8. *If $f \in L^1_{\text{loc}}[0, T) \subset \mathscr{E}$ is a nonpositive distribution. Then, $f \leq 0$ almost everywhere with respect to Lebesgue measure.*

LEMMA 4.9. *Suppose $f_1, f_2 \in \mathscr{G}_c$ for $T \in (0, \infty]$. If $f_1 \leq f_2$ (we mean $f_1 - f_2$ is a nonpositive distribution), and that both $h_1 = J_\gamma f_1$ and $h_2 = J_\gamma f_2$ are functions in $L^1_{\text{loc}}[0, T)$, then*

$$h_1 \leq h_2, a.e.$$

*Proof.* By the definition of $\mathscr{G}_c$, $\mathrm{supp}(f_1) \subset [0,T)$ and $\mathrm{supp}(f_2) \subset [0,T)$.

Suppose the conclusion is not true. Then by Lemma 4.8, there exists $\varphi \in C_c^\infty(-\infty, T)$, $\varphi \geq 0$ such that

$$\int_0^T (h_1 - h_2)\varphi \, dt > 0.$$

Also, we are able to find $\epsilon > 0$ such that $\mathrm{supp}\,\varphi \subset (-\infty, T - \epsilon)$.

This means

$$\langle J_\gamma f_1, \varphi \rangle > \langle J_\gamma f_2, \varphi \rangle.$$

By Proposition 2.9, we can find an extension operator $K_n^T$ so that

$$\langle g_\gamma * \tilde{f}_1, \varphi \rangle > \langle g_\gamma * \tilde{f}_2, \varphi \rangle.$$

where $\tilde{f}_1 = K_n^T f_1$ and $\tilde{f}_2 = K_n^T f_2$, while

$$\langle \tilde{f}_i, \tilde{\varphi} \rangle = \langle f_i, \tilde{\varphi} \rangle, \ i = 1, 2, \ \forall \tilde{\varphi} \in C_c^\infty(-\infty, T), \ \mathrm{supp}\,\tilde{\varphi} \subset (-\infty, T - \epsilon].$$

Let $\{\phi_i\}$ be a partition of unity for $\mathbb{R}$. Then, by Definition 2.1, we have

$$\sum_i \left( \left\langle (\phi_i g_\gamma) * \tilde{f}_1, \varphi \right\rangle - \left\langle (\phi_i g_\gamma) * \tilde{f}_2, \varphi \right\rangle \right) > 0.$$

There are only finitely many terms that are nonzero in this sum. Hence, there must exist $i_0$ such that

$$\left\langle (\phi_{i_0} g_\gamma) * \tilde{f}_1, \varphi \right\rangle - \left\langle (\phi_{i_0} g_\gamma) * \tilde{f}_2, \varphi \right\rangle > 0.$$

Denote $\zeta_{i_0}(t) = (\phi_{i_0} g_\gamma)(-t)$. Then,

$$\langle \tilde{f}_1 - \tilde{f}_2, \zeta_{i_0} * \varphi \rangle > 0.$$

$\zeta_{i_0}$ is a positive integrable function with compact support and $\varphi \geq 0$ is compactly supported smooth function. Then, $\zeta_{i_0} * \varphi \geq 0$ and is $C_c^\infty(-\infty, T)$ with the support in $(-\infty, T - \epsilon]$. It then means

$$\langle f_1 - f_2, \zeta_{i_0} * \varphi \rangle > 0.$$

This is a contradiction since we have assumed $f_1 \leq f_2$. □

Now, we introduce the comparison principle, which is important for a priori energy estimates of time fractional PDEs. This is because the energy $E$ of a time fractional PDE usually satisfies some inequality $D_c^\gamma E \leq f(E)$ and we need to bound $E$.

THEOREM 4.10. *Let $f(t,v)$ be a continuous function, locally Lipschitz in $v$ and $\forall t \geq 0$, $x \leq y$ implies $f(t,x) \leq f(t,y)$. Let $0 < \gamma < 1$. Suppose $v_1(t)$ is continuous satisfying*

$$D_c^\gamma v_1(t) \leq f(t, v_1(t)),$$

*where this inequality means $D_c^\gamma v_1 - f(\cdot, v_1(\cdot))$ is a non-positive distribution. Suppose also that $v_2$ is the continuous solution of the equation*

$$D_c^\gamma v_2(t) = f(t, v_2(t)), \quad v_2(0) \geq v_1(0).$$

Then, $v_1 \leq v_2$ on their common interval of existence.

Correspondingly, if
$$D_c^\gamma v_1(t) \geq f(t, v_1(t)),$$
and $v_2$ solves
$$D_c^\gamma v_2(t) = f(t, v_2(t)), \quad v_2(0) \leq v_1(0).$$
Then, $v_1 \geq v_2$ on their common interval of existence.

*Proof.* Fixing $T_1 \in (0, T_b)$, we show that $v_1(t) \leq v_2(t)$ for any $t \in (0, T_1]$. Then, since $T_1$ is arbitrary, the claim follows. By Theorem 4.4, we can find $\epsilon_0 > 0$ such that the continuous solution of the equation with initial data $v_2(0) + \epsilon$, denoted by $v_2^\epsilon$, exists on $(0, T_1]$ whenever $\epsilon \leq \epsilon_0$.

Define $T^\epsilon = \inf\{t > 0 : v_2^\epsilon(t) \leq v_1(t)\}$. Since both $v_1$ and $v_2^\epsilon$ are continuous and $v_2^\epsilon(0) > v_1(0)$, $T^\epsilon > 0$. We claim that $T^* = T_1$. Otherwise, we have $v_2^\epsilon(T^*) = v_1(T^*)$ and $v_1(t) < v_2(t)$ for $t < T^*$. Note that $f(t, v_1)$ is a continuous function. By using Theorem 3.7 and Lemma 4.9:

$$v_1(T^*) = v_1(0) + \frac{1}{\Gamma(\gamma)} \int_0^{T^*} (T^* - s)^{\gamma - 1} D_c^\gamma v_1(s) \, ds \leq v_1(0) + \frac{1}{\Gamma(\gamma)} \int_0^{T^*} (T^* - s)^{\gamma - 1}$$

$$\times f(s, v_1(s)) ds < v_2(0) + \epsilon + \frac{1}{\Gamma(\gamma)} \int_0^{T^*} (T^* - s)^{\gamma - 1} f(s, v_2(s)) ds = v_2^\epsilon(T^*).$$

(The first integral is understood as $J_\gamma D_c^\gamma v_1$. However, the obtained distribution is a continuous function and Lemma 4.9 guarantees that we can have the first inequality.) This is a contradiction. Hence, $v_1(t) < v_2^\epsilon(t)$ for all $t \in (0, T_1]$. Taking $\epsilon \to 0^+$, using the continuity on initial value yields the claim. Then, by the arbitrariness of $T_1$, the first claim follows.

Similar arguments hold for the second claim, except that we perturb $v_2(0)$ to $v_2(0) - \epsilon$ to construct $v_2^\epsilon$. □

We now show another result that may be useful for time fractional PDEs:

PROPOSITION 4.11. *Suppose $f$ is nondecreasing and Lipschitz on any bounded interval, satisfying $f(0) \geq 0$. Then, the continuous solution of $D_c^\gamma v = f(v), v(0) \geq 0$ is nondecreasing on the interval $[0, T_b)$ where $T_b$ is the blowup time given in Proposition 4.6.*

*Proof.* It is clear that $f(v) \geq 0$ whenever $v \geq 0$. We first show that $v(t) \geq v(0)$ for all $t \in [0, T_b)$.

Let $v^\epsilon$ be the solution with initial data $v(0) + \epsilon > v(0)$. Fix $T_1 \in (0, T_b)$. There exists $\epsilon_0 > 0$ such that $\forall \epsilon \in (0, \epsilon_0)$, $v^\epsilon$ is defined on $(0, T_1]$. Define $T^* = \inf\{t \leq T_1 : v^\epsilon(t) \leq v(0)\}$. $T^* > 0$ because $v^\epsilon(0) > v(0)$. We show that $T^* = T_1$. If this is not true, $v^\epsilon(T^*) = v(0)$ and $v^\epsilon(t) > v(0) \geq 0$ for all $t < T^*$. Applying Theorem 3.7 for $t = T^*$ gives $v^\epsilon(T^*) > v(0)$, a contradiction. Hence, $v^\epsilon \geq v(t)$ for all $t \in (0, T_1]$. Taking $\epsilon \to 0$ yields $v(t) \geq v(0)$ for all $t \in [0, T_1]$, since the solution is continuous on initial data by Theorem 4.4. The arbitrariness of $T_1$ concludes the claim.

Now, consider the function sequence
$$D_c^\gamma v^n = f(v^{n-1}), \quad v^n(0) = v(0), \quad v^0 = v(0) \geq 0.$$

All functions are continuous and defined on $[0, \infty)$. Since $v(t) \geq v^0$ on $[0, T_b)$, then $f(v) \geq f(v^0)$, and it follows that $v \geq v^1$ on $[0, T_b)$ by Theorem 3.7. Doing this iteratively, we find that $v \geq v^n$ for all $n \geq 0$ and all $t \in [0, T_b)$.

Theorem 3.7 shows that

$$v^1(t) = v(0) + \frac{f(v(0))}{\Gamma(1+\gamma)} t^\gamma \geq v^0.$$

By Theorem 3.7 again, it follows that $v_2 \geq v_1$ and hence $v^n \geq v^{n-1}$ for all $n \geq 1$. Therefore, $v^n$ is increasing in $n$ and bounded above by $v$ on $[0, T_b)$. Then, $v^n \to \tilde{v}$ pointwise on $(0, T_b)$ and $\tilde{v}$ is nondecreasing. Taking the limit for

$$v^n(t) = v(0) + \frac{1}{\Gamma(\gamma)} \int_0^t (t-s)^{\gamma-1} f(v^{n-1}(s)) ds,$$

and by monotone convergence theorem, we find that $\tilde{v}$ satisfies

$$\tilde{v}(t) = v(0) + \frac{1}{\Gamma(\gamma)} \int_0^t (t-s)^{\gamma-1} f(\tilde{v}^{n-1}(s)) ds.$$

Since $\tilde{v}$ is bounded by $v$ on any closed subinterval of $[0, T_b)$, $\tilde{v}$ is continuous by this integral equation. Since the solution is unique on $[0, T_b)$, it must be $v$.

This said, now we show that $v^n$ is noncreasing in $t$. This is clear by induction if we note this fact: "If $h(t) \geq 0$ is a nondecreasing locally integrable function, then $g_\gamma * h$ is nondecreasing in $t$", which can be verified by direct computation. $v$, as the the limit of noncreasing functions, is noncreasing. □

Clearly, for functions valued in $\mathbb{R}^m, m \in \mathbb{N}_+$, the fractional derivatives are defined by taking derivatives on each component. The results in Theorem 4.4, Proposition 4.6, Proposition 3.11 can be generalized to $\mathbb{R}^m$ easily and we conclude the following:

PROPOSITION 4.12. *Suppose $E(\cdot) \in C^1(\mathbb{R}^m, \mathbb{R})$ is convex and $\nabla E$ is locally Lipschitz continuous. Let $v \in L^\infty_{\text{loc}}([0, T_b), \mathbb{R}^m)$ solve the fractional gradient flow:*

$$D_c^\gamma v = -\nabla_v E(v) \tag{56}$$

*and $T_b$ is the largest existence time. Then*

$$E(v(t)) \leq E(v(0)), \ t \in [0, T_b).$$

*Similarly, under the same conditions except that the equation is of the form:*

$$D_c^\gamma v = J \nabla_v E(v), \tag{57}$$

*where $J$ is an anti-Hermitian constant operator, we obtain*

$$E(v(t)) \leq E(v(0)), \ t \in [0, T_b).$$

*Proof.* Consider the time fractional gradient flow problem first. Since $\nabla_v E(v)$ is locally Lipschitz, by Theorem 4.4, $v$ is continuous. By the equation, we find that $D_c^\gamma v$ is continuous.

We now take a mollifier $\zeta \in C_c^\infty(\mathbb{R}^m)$ and define

$$E^\epsilon = E * \left( \frac{1}{\epsilon^m} \zeta(\frac{v}{\epsilon}) \right).$$

Let $v^\epsilon$ be the solution to the corresponding fractional ODE with $E$ replaced by $E^\epsilon$. Using the techniques in the proof of Theorem 4.4, $v^\epsilon$ exists on $[0, T]$ when $\epsilon$ is small enough and $v^\epsilon \to v$ in $C([0, T])$ as $\epsilon \to 0$. $v^\epsilon$ satisfies

$$v^\epsilon(t) = v(0) - \frac{1}{\Gamma(\gamma)} \int_0^t (t-s)^{\gamma-1} \nabla_v E^\epsilon(v^\epsilon(s)) ds$$

by Proposition 4.2. Using [22, Thm. 1], we have $v^\epsilon \in C[0,T] \cap C^1(0,T]$ for $\epsilon$ sufficiently small. By Proposition 3.11,

$$D_c^\gamma E(v^\epsilon) \leq \nabla_v E^\epsilon(v^\epsilon) \cdot D_c^\gamma v^\epsilon = -|\nabla_v E^\epsilon(v^\epsilon)|^2.$$

Since $T \in (0, T_b)$ is arbitrary, taking $\epsilon \to 0$, we have in the distributional sense that

$$D_c^\gamma E(v) \leq -|\nabla_v E(v)|^2 \leq 0.$$

By Theorem 4.10 (the function is $f(t, E(v)) = 0$), we conclude that

$$E(v(t)) \leq E(v(0)).$$

For the second case, we have similarly in the distributional sense that

$$D_c^\gamma E(v) \leq \nabla_v E(v) \cdot D_c^\gamma u = \nabla_v E \cdot (J \nabla_v E) = 0.$$

Theorem 4.10 again yields
$$E(v(t)) \leq E(v(0)). \qquad \square$$

As a straightforward corollary, we have

COROLLARY 4.13. *If $\lim_{|v| \to \infty} E(v) = \infty$ and the conditions in Proposition 4.12 hold, then the solutions exist globally. In other words, $T_b = \infty$.*

The fractional Hamiltonian system $D_c^\gamma v = J \nabla_v E(v)$ can be found in [23, 26, 5]. These equations were introduced for systems of nonconservative forces with Lagrangian containing fractional orders. We introduced these systems here only for mathematical study without claiming the physical significance. The fractional Hamiltonian system can be rewritten as $v(t) = v(0) + J_\gamma(J \nabla_v E(v))$ by Theorem 3.7, which is of the Volterra type $v_0 \in v(t) + b * (Av)$. The general Volterra equations with completely positive kernels and $m$-accretive $A$ operators have been discussed in [4], and the solutions have been shown to converge to the equilibrium. In the fractional Hamiltonian system, $-J \nabla_v E(v)$ is not $m$-accretive and it is not clear whether the solutions converge to one equilibrium satisfying $\nabla_v E(v) = 0$ or not. However, for the following simple example, the energy function $E$ indeed dissipates and the solutions converge to the equilibrium (this example is also discussed in [26]):

**Example**: Consider $E(p, q) = \frac{1}{2}(p^2 + q^2)$ and

$$D_c^\gamma q = \frac{\partial E}{\partial p} = p,$$
$$D_c^\gamma p = -\frac{\partial E}{\partial q} = -q,$$

where we assume $\gamma < 1/2$, and the initial conditions are $p(0) = p_0, q(0) = q_0$.

Applying Theorem 4.4 for $v(t) = (p, q)$ and Corollary 4.5, we find that both $p$ and $q$ are continuous functions, and exist globally. By Lemma 3.5, $D_c^\gamma(D_c^\gamma q) = D_c^{2\gamma} q - p_0 g_{1-\gamma}$. Applying $D_c^\gamma$ on the first equation yields

$$D_c^{2\gamma} q - p_0 g_{1-\gamma} = D_c^\gamma p = -q.$$

By Proposition 4.3, we find

$$q(t) = q_0 \beta_{2\gamma}(t) - p_0 g_{1-\gamma} * (\beta'_{2\gamma})(t) = q_0 \beta_{2\gamma}(t) - p_0 D_c^\gamma \beta_{2\gamma}(t),$$

where $\beta_{2\gamma}(t) = E_{2\gamma}(-t^{2\gamma})$ is defined as in Proposition 4.3. Note that $g_{1-\gamma} * (\beta'_{2\gamma}) = D_c^\gamma \beta_{2\gamma}$ is due to Proposition 3.6.

Using the second equation, we find
$$p(t) = p_0 + J_\gamma(-q) = p_0 \beta_{2\gamma}(t) - q_0 J_\gamma \beta_{2\gamma}(t).$$

Actually, from the equation of $p$, $D_c^{2\gamma} p = -p - q(0) g_{1-\gamma}$, we find
$$p(t) = p_0 \beta_{2\gamma}(t) + q_0 D_c^\gamma \beta_{2\gamma}(t).$$

Since $\beta_{2\gamma}$ solves the equation $D_c^{2\gamma} v = -v$, we see $D_c^\gamma \beta_{2\gamma} = -J_\gamma \beta_{2\gamma}$. Those two expressions for $p(t)$ are identical. Hence, we find that

(58) $$E(t) = E(0)(\beta_{2\gamma}^2(t) + (J_\gamma \beta_{2\gamma})^2(t)).$$

(Note that $\beta_{2\gamma}$ and $J_\gamma \beta_{2\gamma}$ are the solutions to the following two equations respectively:
$$D_c^{2\gamma} v = -v, \quad v(0) = 1,$$
$$D_c^{2\gamma} v = -v + g_{1-\gamma}, \quad v(0) = 0.$$

Unlike the corresponding ODE system where the two components are both solutions to $v'' = -v$, here since $D_c^{2\gamma} = -v$ only has one solution for an initial value, the two functions are from different equations.) By the series expression of $E_{2\gamma} = E_{2\gamma,1}$ [16]:

(59)
$$\beta_{2\gamma}(t) = E_{2\gamma}(-t^{2\gamma}) = \sum_{n=0}^{\infty}(-1)^n g_{2n\gamma+1}(t),$$
$$J_\gamma \beta_{2\gamma}(t) = \sum_{n=0}^{\infty}(-1)^n g_{(2n+1)\gamma+1}(t) = t^\gamma E_{2\gamma,\gamma+1}(-t^{2\gamma}).$$

According to the asymptotic behavior listed in [12, eq. (7)],

(60) $$E_{2\gamma,\rho}(-t^{2\gamma}) \sim -\sum_{k=1}^{p} \frac{(-1)^k t^{-2\gamma k}}{\Gamma(\rho - 2\gamma k)}, \quad 0 < \gamma < 1.$$

As examples, $E_{1/2,1}(-t^{1/2}) = e^t \left(1 - \frac{2}{\sqrt{\pi}} \int_0^{t^{1/2}} \exp(-s^2) ds\right) \sim 1/t^{1/2}$. $E_{1,1}(-t) = e^{-t}$ decays exponentially, while in equation (60), $\Gamma(1-k) = \infty$ for all $k = 1, 2, 3, \ldots$. Note that one cannot take the limit $\gamma \to 1$ in Equation (60) and conclude that $E_{2,1}(-t^2)$ and $tE_{2,2}(-t^2)$ both decay exponentially. Actually, $\beta_2 = E_{2,1}(-t^2) = \cos(t)$ and $J_1 \beta_2 = t E_{2,2}(-t^2) = \sin(t)$. The singular limit is due to an exponentially small term $C_1 \exp(C_2(\gamma - 1)t)$ in $E_{2\gamma,\rho}$.

If $\gamma < 1/2$, then $\Gamma(\beta - 2\gamma) \neq \infty$ and the leading order behavior is $t^{-2\gamma}$. Hence, $J_\gamma \beta_{2\gamma}(t) \sim t^{-\gamma}$ and $E(t)$ decays like $t^{-2\gamma}$. Following the same method, one can use the last statement in Lemma 3.5 to solve the case $\gamma = 1/2$ and find that $\gamma = 1/2$ case is still right: $\beta_1 = E_{2,1}(-t) = e^{-t}$ decays exponentially fast while $t^{1/2} E_{2,2}(-t)$ decays like $t^{-1/2}$. Since $E(t)$ decays to zero, the solution must converge to $(0,0)$.

Whether this is true for $1/2 < \gamma < 1$ is interesting, which we leave for future study.

REMARK 7. *According to Theorem 4.4, the system has a unique solution for any $(p(0), q(0)) \in \mathbb{R}^2$ and $0 < \gamma < 1$. We have only considered $0 < \gamma < 1/2$ here just because we can find the solution easily for these cases. If one defines the Caputo derivatives for $1 < \alpha < 2$ ($\alpha = 2\gamma$) consistently, we guess that the expressions above are still correct.*

**Acknowledgements.** L. Li is grateful to Xianghong Chen and Xiaoqian Xu for discussion on Sobolev spaces. The work of J.-G Liu is partially supported by KI-Net NSF RNMS11-07444 and NSF DMS-1514826.